\documentclass[12pt,a4paper,oneside]{article}%
\usepackage{amsfonts}
\usepackage{graphicx}
\usepackage{amsmath,amssymb}
\usepackage{color}
\usepackage{amssymb}
\usepackage{amsmath}
\usepackage[T1]{fontenc}
\usepackage{amsthm}
\usepackage[margin=1in]{geometry}
\usepackage{amssymb}
\usepackage{amsfonts}
\usepackage{graphicx}
\usepackage{amsmath}%
\providecommand{\U}[1]{\protect\rule{.1in}{.1in}}

\DeclareMathOperator{\sign}{sign}

\numberwithin{equation}{section}

\newcommand{\p}{ { \mathbb P} }

\newcommand\E{{\mathbb {E}}}

\newtheorem {Lemma}{Lemma}[section]

\newtheorem {Proposition}{Proposition}[section]
\newtheorem {Corollary}{Corollary}[section]

\newtheorem{Definition}{Definition}[section]

\newtheorem{Remark}{Remark}[section]

\newcommand\F{{\mathcal {F}}}

\newcommand\beq{\begin{equation}}
\newcommand\eeq{\end{equation}}

\def\U{\mathcal{U}}
\def\F{\mathcal{F}}

\begin{document}

\title{Density estimation for $\beta$-dependent sequences}

\date{}

\author{J\'er\^ome Dedecker$^a$ and Florence Merlev\`ede$^b$}
\maketitle

{\small{$^a$ Universit\'e Paris Descartes, Sorbonne Paris Cit\'e,  Laboratoire MAP5 (UMR 8145).

Email: jerome.dedecker@parisdescartes.fr \bigskip

$^b$ Universit\'{e} Paris-Est, LAMA (UMR 8050), UPEM, CNRS, UPEC.

Email: florence.merlevede@u-pem.fr}}

 \bigskip

{\abstract{We study the ${\mathbb L}^p$-integrated risk of some classical estimators of the density, 
when the observations are drawn from a strictly stationary sequence. The results apply to a large class of sequences, 
which can be non-mixing in the sense of Rosenblatt and long-range dependent.
The main probabilistic tool is a new Rosenthal-type inequality for partial sums of $BV$ functions of the variables.  As an application, we give the rates of convergence of 
regular Histograms, when estimating the  invariant density of a class of expanding maps of the unit interval with a neutral fixed point at 
zero. These Histograms are plotted in the section devoted to the simulations.
}}

\bigskip
{\small{

\smallskip

}}

\section{Introduction}

In this paper, we have four goals:

\begin{enumerate}
\item We wish to  extend some of the results of Viennet \cite{V} for stationary $\beta$-mixing sequences 
to the much larger class of $\beta$-dependent sequences, as introduced in 
\cite{DP}, \cite{DP2}. Viennet proved that, if the $\beta$-mixing coefficients $\beta(n)$ of a stationary
sequence $(Y_i)_{i \in {\mathbb Z}}$ are such that 
\begin{equation}\label{ViennetCond}
  \sum_{k=0}^\infty k^{p-2} \beta(k) < \infty \, , \quad \text{for some $p\geq 2$,}
\end{equation}
then the ${\mathbb L}^p$-integrated risk of the usual estimators of the density of   $Y_i$ 
behaves as in the independent and identically distributed  (iid) case   (as described in the paper by Bretagnolle and Huber \cite{BH}). 

For Kernel estimators, we shall obtain a complete extension of Viennet's result (assuming only as an extra hypothesis that the Kernel has bounded variation). 
For projection estimators, the situation
is more delicate, because our dependency coefficients cannot always give a good upper bound
for the variance of the estimator (this was already pointed out in \cite{DP}).
However, for estimators based on piecewise polynomials  (including Histograms), the result of Viennet 
can again be fully extended.

\item We shall consider the ${\mathbb L}^p$-integrated risk for any $p\in [1, \infty)$, and not 
only for $p\geq 2$ (which was the range considered in \cite{V}). Two main reasons 
for this: first the case $p=1$ is of particular interest, because it gives some information 
on the total variation between the (possibly signed) measure with density $f_n$ (the estimated
density) and the distribution of $Y_i$. The variation distance is a true distance between measures,
contrary to the ${\mathbb L}^p$-distance between densities, which depends on the dominating measure. Secondly, we have in mind applications to some  classes of dynamical systems
(see point 4 below),
 for
which it is known either that the density has bounded variation over $[0,1]$ or that it is 
non-decreasing on $(0,1]$ (and blows up as $x \rightarrow 0$). In such cases, it turns 
out that the bias of our estimators is well controlled in ${\mathbb L}^1([0,1], dx)$. 

\item We want to know what happens if \eqref{ViennetCond} is not satisfied, or
if $\sum_{k \geq 0} \beta(k) = \infty$ in the case where $p \in [1,2]$. Such results 
are not given in the paper \cite{V}, although Viennet could have done it by 
refining some computations. For $\beta$-dependent sequences and  $p=2$, the situation 
is clear (see \cite{DP}): the rate of convergence of the estimator depends 
on the regularity of $f$ and of the behavior of 
$\sum_{k = 0}^n \beta_{1, Y}(k) $ (the coefficients $\beta_{1,Y}(k)$ will be defined in the next
section, and are weaker than the corresponding $\beta$-mixing coefficients). 
Hence, in that case, the consistency holds as soon as $\beta_{1,Y}(n)$ tends to zero
as $n$ tends to infinity, and one can compute the rates of convergence as soon as 
one knows the asymptotic behavior of $\beta_{1,Y}(n)$. This is the kind of result
we want to extend to any $p \in [1, \infty)$. Once again, we have precise motivations 
for this, coming from dynamical systems that can exhibit long-range dependence (see point 4 below).

\item As already mentioned, our first motivation was to study the robustness of 
the usual estimators of the density, showing that they apply to a  larger class of dependent processes
than in \cite{V}. But our second main objective was to be able to visualize 
the invariant density of the iterates of expanding maps of the unit interval. For uniformly expanding
maps the invariant density has bounded variation, and one can estimate it  in 
${\mathbb L}^1([0,1], dx)$
at the usual rate
$n^{-1/3}$ by using an appropriate Histogram (see Subsection \ref{UE}). The case of the intermittent 
map $T_\gamma$ 
(as defined in \eqref{LSV} for $\gamma \in (0,1)$) is even more interesting. In that case, one knows 
that the invariant density $h_\gamma$  is equivalent 
to the density $x\rightarrow (1-\gamma) x^{-\gamma}$ on $(0,1)$ (see the inequality \eqref{equiv}). 
Such a map exhibits long-range dependence as soon as $\gamma \in (1/2,1)$, but we can 
use our upper bound for the random part  + bias of a regular Histogram, to compute 
the appropriate number of breaks of the Histogram (more precisely we give the order of the  number of breaks
as a function of  $n$, up to an unknown constant). In Figure 5 of Section \ref{Sec:Simu}, we plot
the Histograms of the invariant density  $h_\gamma$, when $\gamma=1/4$ (short-range dependent
case),  $\gamma=1/2$ (the boundary case),  and $\gamma=3/4$ (long-range dependent case).
\end{enumerate}

One word about the main probabilistic tool. As shown in \cite{V}, to control the
${\mathbb L}^p$-integrated risks for $p>2$, the appropriate tool is a precise 
Rosenthal-type inequality. Such an inequality is not easy to prove in the 
$\beta$-mixing case, and the $\beta$-dependent case is even harder to handle
because we cannot use Berbee's coupling (see \cite{B}) as in \cite{V}. 
A major step to get a good Rosenthal bound has been made by Merlev\`ede
and Peligrad \cite{MP}: they proved a very general 
inequality involving only conditional expectations of the random variables 
with respect to the past $\sigma$-algebra, which can be applied to many 
situations. However, it does not fit completely to our context, and leads to small
losses when applied to kernel estimators (see Section 5 in \cite{MP}). In Section \ref{Rosineq}, we shall prove a taylor-made inequality, in the spirit 
of that of Viennet but expressed in terms of our weaker coefficients. This inequality
will give the complete extension of Viennet's results for Kernel estimators and 
estimators based on piecewise polynomials, when $p>2$.

\section{A Rosenthal inequality for $\beta$-dependent sequences}\label{Rosineq}

From now, $(Y_i)_{i \in {\mathbb Z}}$ is a strictly stationary sequence of real-valued random variables.
We define the $\beta$-dependence coefficients of $(Y_i)_{i \in {\mathbb Z}}$  as in \cite{DP2}:

\begin{Definition}\label{beta}
Let $P$ be the law of $Y_0$ and $P_{(Y_i, Y_j)}$ be the law of $(Y_i, Y_j)$. Let ${\mathcal F}_{\ell} = \sigma( Y_i, i \leq \ell)$, let
 $P_{Y_k| {\mathcal
F}_\ell}$ be the conditional distribution of  $Y_k$ given ${\mathcal F}_\ell$, and  let $P_{(Y_i,
Y_j)|{\mathcal F}_\ell}$ be the conditional distribution of $(Y_i,Y_j)$ given ${\mathcal
F}_\ell$. Define the functions  $$f_t={\bf 1}_{]-\infty, t]}\, \quad \text{and} \quad f_t^{(0)}=f_t-P(f_t) \, , $$ and the
 the random variables
\begin{eqnarray*}
b_{ \ell} (k) &=& \sup_{t \in {\mathbb R}}\left |{P}_{Y_{k}|{\mathcal F}_\ell}(f_t)-P(f_t) \right |\, , \\
b_{ \ell} ( i,j)&=& \sup_{(s,t) \in {\mathbb R}^2} \left |P_{(Y_i, Y_j)|{\mathcal F}_\ell}\left (f_t^{(0)}\otimes f_s^{(0)} \right )-P_{(Y_i, Y_j)}\left (f_t^{(0)}\otimes f_s^{(0)}\right ) \right |\, .
\end{eqnarray*}
Define now the coefficients
\[
\beta_{1,Y}(k)= {\mathbb E}(b_{0}( k))\,  \text{ and } \,  \beta_{2,Y}(k) =\max \left \{\beta_1(k), \sup_{i>j\geq k} {\mathbb E}((b_0( i,j))) \right \}\, .
\]
\end{Definition}
These coefficients are weaker than the usual $\beta$-mixing coefficients of $(Y_i)_{i \in {\mathbb Z}}$. Many examples of non-mixing process 
for which $\beta_{2,Y}(k)$ can be computed are given in \cite{DP2}. 
Some of these examples will be studied in Sections \ref{Sec:EM} and \ref{Sec:Simu}.

Let us now give the main probabilistic tool of the paper. It is a Rosenthal-type inequality for partial sums of $BV$ functions of $Y_i$ (as usual, $BV$ means ``of bounded variation''). We shall use it to control 
the random part of the ${\mathbb L}^p$ integrated risk of the estimators of the density of $Y_i$ when $p>2$. The proof of this inequality  is given in Subsection \ref{proofrose}. Is is quite delicate, and
relies  on two intermediate results (see Subsection \ref{devineq}). 

In all the paper, we shall use  the  notation
 $a_n \ll b_n$, which means that there exists a positive  constant $C$ not
depending on $n$ such that  $a_n \leq  Cb_n $, for all positive integers  $n$.

\begin{Proposition} \label{corRosenthalbeta} Let
$(Y_i)_{i \in {\mathbb Z}}$ be a strictly stationary sequence of real-valued random variables.
For any $p > 2$ and any positive integer $n$, there exists a non-negative ${\mathcal F}_0$-measurable  random variable $A_0(n,p)$
satisfying $\E(A_0(n,p)) \leq  \sum_{k=1}^{n} k^{p-2} \beta_{2,Y} (k)$ and such that: for any $BV$ function $h$ from ${\mathbb R}$
to ${\mathbb R}$, letting 
 $X_i=h(Y_i) - {\mathbb E} (h(Y_i))$ and $S_n = \sum_{k=1}^n X_k$, we have
  \begin{multline}
  \label{ineRosbeta}
  \E \left( \sup_{ 1 \leq k \leq n}  |
  S_k |^p \right)  \ll  n^{p/2} \left ( \sum_{i=0}^{n-1} | {\rm Cov} (X_0, X_i)| \right )^{p/2}  
 +  n \Vert dh \Vert^{p-1} \E \left (   |  h(Y_0) |  A_0(n,p)  \right ) \\
+ n  \Vert dh \Vert^{p-1} \E \left (|h(Y_0)|  \right )  \sum_{k=0}^{n} (k+1)^{p-2} \beta_{2,Y} (k) \, ,
  \end{multline}
where $\Vert dh \Vert$ is the total variation norm of the measure $dh$.
\end{Proposition}

\begin{Remark}
For $p=2$, using Proposition 1 in \cite{DR}, the inequality can be simplified as follows.
There exists a non-negative ${\mathcal F}_0$-measurable  random variable $A_0(n,2)$
satisfying $\E(A_0(n,2)) \leq  \sum_{k=1}^{n}  \beta_{1,Y} (k)$  and such that: for any $BV$ function $h$ from ${\mathbb R}$
to ${\mathbb R}$, letting 
 $X_i=h(Y_i) - {\mathbb E} (h(Y_i))$ and $S_n = \sum_{k=1}^n X_k$, we have
  \begin{equation}
  \label{ineRosbeta2}
  \E \left( \sup_{ 1 \leq k \leq n}  
  S_k ^2 \right)  \ll  
  n \Vert dh \Vert
 \E \left (   |  h(Y_0) |  (1+ A_0(n,2)) \right ) \, .
  \end{equation}
\end{Remark}

\section{${\mathbb L}^p$-integrated risk for Kernel estimators}\label{mainsec}

 Let $(Y_{i})_{i \in \mathbb Z}$ be a stationary sequence
with unknown marginal density $f$. In this section, we wish to build an estimator of $f$ based on the variables $Y_1, \ldots, Y_n$.

Let $K$ be a bounded-variation function in ${\mathbb L}^1({\mathbb R}, \lambda)$, where $\lambda$ is the 
Lebesgue measure. Let $\|d K \|$ be the 
variation norm of the measure $d K$, and $\| K \|_{1, \lambda}$ be the ${\mathbb L}^1$-norm of 
$K$ with respect to $\lambda$. 

Define then
\[
X_{k,n}(x)=K(h_{n}^{-1}(x-Y_{k})) \quad 
\text{and}
\quad f_{n}(x)=\frac{1}{nh_{n}}
\sum_{k=1}^{n}X_{k,n}(x)\,,
\]
where $(h_{n})_{n\geq1}$ is a sequence of positive real numbers.

The following proposition gives an upper bound of the term
\begin{equation}\label{random}
{\mathbb{E}} \left (\int_{\mathbb{R}}|f_{n}(x)-{\mathbb{E}}(f_{n}(x))|^{p}dx \right ) 
\end{equation}
when 
$p > 2$ and $f \in {\mathbb L}^{p}({\mathbb R}, \lambda)$.

\begin{Proposition}\label{p>2}
Let $p > 2$, and assume that $f$ belongs to 
${\mathbb L}^{p}({\mathbb R}, \lambda)$. 
Let 
\begin{equation}\label{betaquant}
V_{1, p, Y}(n)=\sum_{k=0}^n (k+1)^{p-2} \beta_{1, Y}(k) 
\quad 
\text{and}
\quad 
V_{2, p,  Y}(n)=\sum_{k=0}^n (k+1)^{p-2} \beta_{2, Y}(k) 
\, .
\end{equation}
The following upper bounds holds
\begin{multline}\label{mainineq}
{\mathbb{E}} \left (\int_{\mathbb{R}}|f_{n}(x)-{\mathbb{E}}(f_{n}(x))|^{p}dx \right )\ll
\left (\frac{\|dK\|\|K\|_{1, \lambda}
\|f\|_{p, \lambda}^{\frac{(p-2)}{(p-1)}}
\left( V_{1, p, Y}(n) \right )^{\frac{1}{(p-1)}}}{n h_n} \right )^{\frac p 2}
\\ 
+ \frac{1}{(nh_{n})^{p-1}} \|dK\|^{p-1} \|K\|_{1, \lambda} 
V_{2, p, Y}(n)  \, .
\end{multline}
\end{Proposition}

\begin{Remark}
Note  that, if $n h_n \rightarrow \infty$ as $n \rightarrow \infty$ and 
$\sum_{k=0}^\infty (k+1)^{p-2} \beta_{2, Y}(k)< \infty$,
 then it follows from Proposition \ref{p>2} that
\begin{equation} \label{viennet}
\limsup_{n \rightarrow \infty} \
(nh_n)^{p/2} 
{\mathbb{E}} \left (\int_{\mathbb{R}}|f_{n}(x)-{\mathbb{E}}(f_{n}(x))|^{p}dx \right )\leq C \|dK\|^{\frac p 2} \|K\|_{1, \lambda}^{\frac p 2}
\|f\|_{p, \lambda}^{\frac{p(p-2)}{2(p-1)}} \, ,
\end{equation}
for some positive constant $C$.
Note that \eqref{viennet} is comparable to the upper bound 
obtained by Viennet \cite{V}, with two differences: 
firstly our condition is written in terms of the coefficients $\beta_{2, Y}(n)$ (while Viennet used 
the usual $\beta$-mixing coefficients), and secondly we only require that 
$f$ belongs to 
${\mathbb L}^{p}({\mathbb R}, \lambda)$ (while Viennet
assumed that $f$ is bounded). When $p\geq 4$, an upper bound similar to \eqref{viennet} is given 
in \cite{MP}, Proposition 33 Item (2), under the slightly stronger condition
$\beta_{2, Y}(n)=O(n^{-(p-1 + \varepsilon)})$ for some $\varepsilon>0$.
\end{Remark}

\begin{Remark}
The first term in the upper bound of Proposition \ref{p>2} has been obtained by assuming only that 
$f$ belongs to ${\mathbb L}^{p}({\mathbb R}, \lambda)$. As will be clear from the proof, a better upper bound can 
be obtained by assuming that $f$ belongs to ${\mathbb L}^{q}({\mathbb R}, \lambda)$ for $q>p$. For instance, if 
$f$ is bounded, the first term of the upper bound can be replaced by 
$$
 \|f\|_{\infty}^{\frac p 2 -1} 
\left (\frac{\|dK\|\|K\|_{1, \lambda}}{n h_n}
 \right )^{\frac p 2} \sum_{k=0}^n (k+1)^{\frac p 2-1} \beta_{1, Y}(k) \, .
$$
This can lead to a substantial improvement of the upper bound of \eqref{random}, for instance in the cases where
$\sum_{k=0}^\infty (k+1)^{p-2} \beta_{2, Y}(k)= \infty$ but $\sum_{k=0}^\infty (k+1)^{p/2-1} \beta_{1, Y}(k)< \infty$.
\end{Remark}

We now give an upper bound of the same quantity when $1 \leq p \leq 2$. Note that the case $p=1$ is of special interest, 
since it enables to get the rate of 
convergence to the unknow probability $\mu$ (with
density $f$) for the total variation distance. 

\begin{Proposition}\label{p=1}
As in \eqref{betaquant}, let  $V_{1, 2,  Y}(n)=\sum_{k=0}^n  \beta_{1, Y}(k)$.  
The following upper bounds hold
\begin{enumerate}
\item For $p=2$,
$\displaystyle
{\mathbb{E}} \left (\int_{\mathbb{R}}|f_{n}(x)-{\mathbb{E}}(f_{n}(x))|^2 dx \right )\ll
\frac{1}{nh_{n}} \|dK\|
\|K\|_{1, \lambda} 
 V_{1, 2, Y}(n) \, .
$
\item Let $ 1 \leq p <2$ and
$\alpha>1, q>1$. Let also
$
U_{1, q, Y}(n)=\sum_{k=0}^n (k+1)^{\frac 1{q-1}} \beta_{1, Y}(k) 
$.  If
$$
 M_{\alpha q, p}(f) :=\int  |x|^{\frac{\alpha q(2-p)}{p}} f(x) dx < \infty \quad \text{and}
 \quad 
 M_{\alpha, p}(K):= \int |x|^{\frac{\alpha(2-p)}{p}}|K(x)| dx < \infty \, ,
$$
then
\begin{multline*}
{\mathbb{E}} \left (\int_{\mathbb{R}}|f_{n}(x)-{\mathbb{E}}(f_{n}(x))|^p dx \right )\ll
 \left (\frac{ \|dK\| \|K\|_{1, \lambda} 
(M_{\alpha q, p}(f))^{\frac 1 q} 
(U_{1, q,  Y}(n))^{\frac{q-1}{q}}}{ nh_{n} }\right )^{\frac p 2}
\\
+  \left(\frac{\|dK\|
\left (\|K\|_{1, \lambda} + h_n^{\frac{\alpha (2-p)}{p}} M_{\alpha, p}(K)\right )
 V_{1, 2, Y}(n)}{nh_{n}} \right)^{\frac p 2} \, .
\end{multline*}
\end{enumerate}
\end{Proposition}

\begin{Remark}
Note first that Item 1 of Proposition \ref{p=1} is due to Dedecker and Prieur \cite{DP}.
For $ p  \in [1,2)$, 
It follows from Item 2   that 
if ${\mathbb E}(|Y_0|^{\alpha q(2-p)/p})< \infty$ for $\alpha >1, q>1$
and if $
\sum_{k=0}^\infty (k+1)^{1/(q-1)} \beta_{1, Y}(k)
< \infty$, then 
$$
\limsup_{n \rightarrow \infty} (n h_n)^{\frac p 2}
{\mathbb{E}} \left (\int_{\mathbb{R}}|f_{n}(x)-{\mathbb{E}}(f_{n}(x))|^p dx \right ) \leq C \|dK\|^{\frac p 2} \|K\|_{1, \lambda}^{\frac p 2} \left( 1+
(M_{\alpha q,p}(f))^{\frac {p} {2q}}\right) 
$$
for some positive constant $C$, 
provided $nh_n \rightarrow \infty$ and $h_n \rightarrow 0$ as $n \rightarrow \infty$.
As will be clear from the proof, this upper bound 
remains true when $q=\infty$, that is when 
$\|Y_0\|_\infty < \infty$ and 
$\sum _{k=0}^ \infty \beta_{1, Y}(k)
< \infty$.
\end{Remark}

With these two propositions, one can get the rates of convergence of $f_n$ to $f$,  when $f$ belongs
to the generalized Lipschitz spaces Lip$^*(s, {\mathbb L}^p({\mathbb
R}, \lambda))$ with  $s>0$, as defined in \cite{DVL}, Chapter 2, Paragraph 9. Recall that 
Lip$^*(s, {\mathbb L}^p({\mathbb
R}, \lambda))$
 is a particular case of Besov spaces (precisely
Lip$^*(s, {\mathbb L}^p({\mathbb
R}, \lambda)) = B_{s,p, \infty}({\mathbb R})$). Moreover, 
if $s$ is a positive integer, 
Lip$^*(s, {\mathbb L}^p({\mathbb
R}, \lambda))$ contains the Sobolev space 
$W^s({\mathbb L}^p({\mathbb
R}, \lambda))$ if $p>1$ and $W^{s-1}(BV)$ if 
$p=1$ (see again \cite{DVL}, 
Chapter 2, paragraph 9). 
Recall that, if $s$  is a positive integer the space 
$W^s({\mathbb L}^p({\mathbb
R}, \lambda))$ (resp. $W^s(BV)$)  is the space of functions 
for which $f^{(s-1)}$ is absolutely continuous, with almost 
everywhere derivative $f^{(s)}$ belonging to 
${\mathbb L}^p({\mathbb R}, \lambda)$ (resp. $f^{(s)}$
has bounded variation).

Let $K_h(\cdot)=h^{-1}K(\cdot/h)$, and $r$ be a
positive integer,  and assume that, 
for any $g$ in $W^r({\mathbb L}^p({\mathbb R}, \lambda
))$,
\begin{equation}\label{bre}
 \int_{\mathbb R} |g(x)-g * K_h (x) |^p dx \leq C_1
 h^{pr} \|g^{(r)}\|_{p, \lambda}^p \, , 
\end{equation}
for some constant $C_1$ depending only on 
$r$. For instance, \eqref{bre} is 
satisfied for any Parzen kernel of order $r$ (see
Section 4 in \cite{BH}).

From \eqref{bre} and Theorem 5.2 in \cite{DVL}, we infer that, for any 
$g \in {\mathbb L}^p({\mathbb R}, \lambda)$,
$$
 \int_{\mathbb R} |g(x)-g * K_h (x) |^p dx \leq C_2
 (\omega_r(g, h)_p)^p \, , 
$$
for some constant $C_2$ depending only on
$r$, where $\omega_r(g, \cdot)_p$ is the $r$-th modulus of
regularity of $g$ in $ {\mathbb L}^p({\mathbb R}, \lambda)$ as defined in \cite{DVL}, Chapter 2, Paragraph 7.
This last inequality implies that, if $g$
belongs to Lip$^*(s, {\mathbb L}^p({\mathbb
R}, \lambda))$ for any $s \in [r-1, r)$, then
\begin{equation}\label{bias}
 \int_{\mathbb R} |g(x)-g * K_h (x) |^p dx \leq C_2
 h^{ps} \|g\|_{\mathrm{Lip}^*(s, {\mathbb L}^p({\mathbb
R}, \lambda))}^p \, .
\end{equation}

Combining Proposition \ref{p>2}  or 
\ref{p=1} with the control of 
the bias given in \eqref{bias}, we obtain the following
upper bounds for the ${\mathbb L}^p$-integrated
risk of the kernel estimator. 

Let  $K$ be  a bounded variation 
function in ${\mathbb L}^1({\mathbb
R}, \lambda)$, and assume that $K$ satisfies 
\eqref{bre} for some positive integer $r$.

\begin{itemize}
\item Let $p \geq  2$ and assume that 
$\sum_{k=0}^n (k+1)^{p-2} \beta_{2, Y}(k)=O(n^{\delta(p-1)})$ for some 
$\delta \in [0,1)$. 
Assume that $f$ belongs to Lip$^*(s, {\mathbb L}^p({\mathbb
R}, \lambda))$ for $s \in [r-1,r)$  or to $W^s({\mathbb L}^p({\mathbb
R}, \lambda))$ for $s=r$. Then, taking 
$h_n=Cn^{-(1-\delta)/(2s +1)}$, 
$$
{\mathbb{E}} \left (\int_{\mathbb{R}}|f_{n}(x)-f(x)|^{p}dx \right )=O\left (n^{\frac{-ps(1-\delta)}{2s +1}}
\right ) \, .
$$
Hence, if $\sum_{k=0}^\infty (k+1)^{p-2} \beta_{2, Y}(k)<\infty$ (case $\delta=0$),  we obtain the same rate as in the iid situation. This result generalizes the 
result of Viennet \cite{V}, who obtained the same rates
under the condition $\sum_{k=0}^\infty (k+1)^{p-2} \beta (k)<\infty$, where the $\beta(k)$'s are the usual 
$\beta$-mixing coefficients.

\item 
Let $1\leq p <2$.
Assume that $Y_0$ has a moment of order $q(p-2)/p + \varepsilon$ for some $q>1$ and $\varepsilon>0$, and that $
\sum_{k=0}^n (k+1)^{ 1/(q-1)} \beta_{1, Y}(k)
=O(n^{\delta q/(q-1)})$ for some 
$\delta \in [0,1)$. 
Assume that $f$ belongs to Lip$^*(s, {\mathbb L}^p({\mathbb
R}, \lambda))$ for $s \in [r-1,r)$  or to $W^s({\mathbb L}^p({\mathbb
R}, \lambda))$ for $s=r$. Then, taking 
$h_n=Cn^{-(1-\delta)/(2s +1)}$, 
$$
{\mathbb{E}} \left (\int_{\mathbb{R}}|f_{n}(x)-f(x)|^pdx \right )=O\left (n^{\frac{-ps(1-\delta)}{2s +1}}
\right ) \, .
$$
Hence, if $\sum_{k=0}^\infty (k+1)^{1/(q-1)}\beta_{1, Y}(k)<\infty$ (case $\delta=0$),  we obtain the same rate as in the iid situation.

Let us consider the particular case where  $p=1$.  Let $\mu$ be the  probability measure with density $f$,  and
let $\hat \mu_n$ be the  random (signed) measure with density
$f_n$. We have just proved that
$$
{\mathbb E} \|\hat \mu_n -\mu\|=O\left (n^{\frac{-s(1-\delta)}{2s +1}}
\right )
$$
(recall that $\|\cdot \|$ is the variation norm). It is 
an easy exercice to modify $\hat \mu_n$ in order to 
get a random probability measure $\mu_n^*$ that converges to $\mu$  at the same rate (take the 
positive part of $f_n$ and renormalize). 
\end{itemize}

\section{${\mathbb L}^p$-integrated risk
for  estimators based on piecewise polynomials} \label{mainsec2}

 Let $(Y_{i})_{i \in \mathbb Z}$ be a stationary sequence
with unknown marginal density $f$. In this section, we wish to estimate  $f$ on a compact interval $I$
with the help of  the variables $Y_1, \ldots, Y_n$.
Without loss of generality, we shall assume here that
$I=[0,1]$. 

We shall consider the piecewise polynomial basis
on a a regular partition of $[0,1]$, defined as follows. 
Let $(Q_i)_{1 \leq i \leq r+1}$ be an orthonormal
basis of the space of polynomials of order $r$ on 
$[0,1]$, and define the function $R_i$ on 
${\mathbb R}$ by: $R_i(x)=Q_i(x)$ if $x$ belongs to
$]0,1]$ and $0$ otherwise. Consider now the regular partition of $]0,1]$ into $m_n$ intervals $(](j-1)/m_n, j/m_n])_{1 \leq j \leq m_n}$. Define the functions 
$\varphi_{i,j}(x)=\sqrt{ m_n} R_i(m_n x-(j-1))$. Clearly the 
family $(\varphi_{i,j})_{1 \leq i \leq r+1}$ is an orthonormal basis of the space of polynomials of
order $r$ on the interval $[(j-1)/m_n, j/m_n]$.  Since the 
supports of $\varphi_{i,j}$ and $\varphi_{k,\ell}$ are 
disjoints for $\ell \neq j$, the family
$(\varphi_{i,j})_{1 \leq i \leq r+1, 1 \leq j \leq m}$ is then 
an orthonormal system of ${\mathbb L}^2([0,1], \lambda)$. The case of regular Histograms corresponds to
$r=0$.

Define then 
$$
X_{i,j,n}= \frac 1 n \sum_{k=1}^n \varphi_{i,j}(Y_k) 
\quad  \text{and} \quad f_n= \sum_{i=1}^{r+1} \sum_{j=1}^{m_n} X_{i,j,n} \varphi_{i,j} \, .
$$

The following proposition gives an upper bound of 
$$
{\mathbb{E}} \left (\int_0^1|f_{n}(x)-{\mathbb{E}}(f_{n}(x))|^{p}dx \right ) 
$$
when 
$p > 2$ and $f{\bf 1}_{[0,1]} \in {\mathbb L}^{p}([0,1], \lambda)$.

\begin{Proposition}\label{p>2bis}
Let $p > 2$, and assume that $f{\bf 1}_{[0,1]}$ belongs to 
${\mathbb L}^{p}([0,1], \lambda)$. 
Let 
$$
C_{1,p}= \sum_{i=1}^{r+1} \|R_i\|_\infty^{\frac{3p}{2}}
\|d R_i \|^{\frac{p}{2}} \quad \text{and} \quad 
C_{2,p}= \sum_{i=1}^{r+1} \|R_i\|_\infty^{p+1}
\|d R_i \|^{p-1} \, .
$$
and recall that 
$
V_{1, p, Y}(n)$ and
$
V_{2, p,  Y}(n)
$
have been defined in \eqref{betaquant}. 
The following upper bounds holds
\begin{multline*}
{\mathbb{E}} \left (\int_0^1|f_{n}(x)-{\mathbb{E}}(f_{n}(x))|^{p}dx \right )\ll
  \left (\frac {m_n} {n} \right )^{\frac p 2} C_{1,p}\left (
\|f{\bf 1}_{[0,1]}\|_{p, \lambda}^{\frac{(p-2)}{(p-1)}}
\left( V_{1, p, Y}(n) \right )^{\frac{1}{(p-1)}}\right )^{\frac p 2}
\\ 
+  \left ( \frac {m_n} {n} \right )^{p-1} 
C_{2,p} V_{2, p, Y}(n)  \, .
\end{multline*}
\end{Proposition}

\begin{Remark}
Note  that, if $n /m_n \rightarrow \infty$ as $n \rightarrow \infty$ and 
$\sum_{k=0}^\infty (k+1)^{p-2} \beta_{2, Y}(k)< \infty$,
 then it follows from Proposition \ref{p>2bis} that
\begin{equation} \label{viennetbis}
\limsup_{n \rightarrow \infty} \
\left ( \frac {n} {m_n}  \right )^{\frac p 2} 
{\mathbb{E}} \left (\int_0^1|f_{n}(x)-
{\mathbb{E}}(f_{n}(x))|^p dx \right )\leq C 
\cdot C_{1,p} 
\|f{\bf 1}_{[0,1]}\|_{p, \lambda}^{\frac{p(p-2)}{2(p-1)}} \, ,
\end{equation}
for some positive constant $C$.
This bound is comparable to the upper bound obtained by Viennet (\cite{V}, Theorem 3.2) 
for the usual $\beta$-mixing coefficients. Note however that Viennet's results is valid
for a much broader class of projection estimators. As a comparison, it seems very difficult to
deal with the trigonometric basis in our setting.
\end{Remark}

We now give an upper bound of the same quantity when $ 1 \leq p \leq 2$.
\begin{Proposition}\label{p=1bis}
Let $ 1 \leq p \leq 2$.
The following upper bounds holds
\begin{equation*}
{\mathbb{E}} \left (\int_0^1|f_{n}(x)-{\mathbb{E}}(f_{n}(x))|^p dx \right )\ll
 \left( \frac  {m_n} {n} \right )^{\frac p 2} 
C_{1,2}^{\frac p 2} \left ( V_{1, 2, Y}(n) \right)^{\frac p 2} \, .
\end{equation*}
\end{Proposition}

\begin{Remark}
Note first that the upper bounf for $p=2$ is due to Dedecker and Prieur \cite{DP}. Note also that,
if $n /m_n \rightarrow \infty$ as $n \rightarrow \infty$ and 
$\sum_{k=0}^\infty  \beta_{1, Y}(k)< \infty$,
 then it follows from Proposition \ref{p=1bis} that
\begin{equation} \label{viennet3}
\limsup_{n \rightarrow \infty} \
\left ( \frac {n} {m_n}  \right )^{\frac p 2} 
{\mathbb{E}} \left (\int_0^1|f_{n}(x)-
{\mathbb{E}}(f_{n}(x))|^p dx \right )\leq C \cdot
C_{1,2}^{\frac p 2}
 \, ,
\end{equation}
for some positive constant $C$.
\end{Remark}

With these two propositions, one can get the  rates of convergence  of $f_n$ to  $f{\bf 1}_{[0,1]}$ when 
$f{\bf 1}_{[0,1]}$ belongs to 
to the generalized Lipschitz spaces Lip$^*(s, {\mathbb L}^p([0,1], \lambda))$ with  $s>0$. 

Applying  the Bramble-Hilbert lemma (see \cite{BrHi}), we know that, for any $f$ such that $ f{\bf 1}_{[0,1]}$ belongs to  $W^{r+1}({\mathbb L}^p([0,1], \lambda))$
$$
 \int_0^1 |f(x)- {\mathbb E}(f_n(x))|^p dx \leq C_1 m_n^{-p(r+1)} \left  \|f^{(r)}{\bf 1}_{[0,1]} \right \|_{p, \lambda} \, ,
$$
for some constant $C_1$ depending only on $r$. From \cite{DVL}, page 359, 
 we know that,  if  $f{\bf 1}_{[0,1]}$ belongs to Lip$^*(s, {\mathbb L}^p([0,1], \lambda))$ and 
 if the degree $r$ is such that $r>s-1$, 
\begin{equation}\label{biasbis}
 \int_0^1 |f(x)- {\mathbb E}(f_n(x))|^p dx \leq C_2 m_n^{-ps} \| f{\bf 1}_{[0,1]}\|_{\mathrm{Lip}^*(s, {\mathbb L}^p([0,1], \lambda))}\, ,
\end{equation}
for some constant $C_2$ depending only on $r$.
Combining Proposition \ref{p>2bis}  or 
\ref{p=1bis} with the control of 
the bias given in \eqref{biasbis}, we obtain the following
upper bounds for the ${\mathbb L}^p$-integrated
risk.

\begin{itemize}
\item Let $p >  2$ and assume that 
$\sum_{k=0}^n (k+1)^{p-2} \beta_{2, Y}(k)=O(n^{\delta(p-1)})$ for some 
$\delta \in [0,1)$. 
Assume that $f{\bf 1}_{[0,1]}$ belongs to Lip$^*(s, {\mathbb L}^p([0,1], \lambda))$ for $s < r+1$ or to  $W^s({\mathbb L}^p([0,1]))$ for $s=r+1$. Then, taking 
$m_n=[C n^{(1-\delta)/(2s +1)}]$, 
$$
{\mathbb{E}} \left (\int_0^1|f_{n}(x)-f(x)|^{p}dx \right )=O\left (n^{\frac{-ps(1-\delta)}{2s +1}}
\right ) \, .
$$
Hence, if $\sum_{k=0}^\infty (k+1)^{p-2} \beta_{2, Y}(k)<\infty$ (case $\delta=0$),  we obtain the same rate as in the iid situation. This result generalizes the 
result of Viennet \cite{V}, who obtained the same rates
under the condition $\sum_{k=0}^\infty (k+1)^{p-2} \beta (k)<\infty$, where the $\beta(k)$'s are the usual 
$\beta$-mixing coefficients.

\item 
Let $1\leq p \leq 2$ and assume that $\sum_{k=0}^n  \beta_{1, Y}(k)=O(n^{\delta})$  for some 
$\delta \in [0,1)$. 
Assume that $f{\bf 1}_{[0,1]}$ belongs to Lip$^*(s, {\mathbb L}^p([0,1], \lambda))$ for $s < r+1$ or to  $W^s({\mathbb L}^p([0,1], \lambda))$ for $s=r+1$. Then, taking 
$m_n=[Cn^{(1-\delta)/(2s +1)}]$, 
$$
{\mathbb{E}} \left (\int_0^1 |f_{n}(x)-f(x)|^pdx \right )=O\left (n^{\frac{-ps(1-\delta)}{2s +1}}
\right ) \, .
$$
Hence, if $\sum_{k=0}^\infty \beta_{1, Y}(k)<\infty$ (case $\delta=0$),  we obtain the same rate as in the iid situation.

Let us consider the particular case where  $p=1$ and $f$ is supported on $[0,1]$.  Let $\mu$ be the  probability measure with density $f$,  and
let $\hat \mu_n$ be the  random   measure with density
$f_n$. We have just proved that
$$
{\mathbb E} \|\hat \mu_n -\mu\|=O\left (n^{\frac{-s(1-\delta)}{2s +1}}
\right )
$$
(recall that $\|\cdot \|$ is the variation norm). 
\end{itemize}

\section{Application to density estimation of expanding maps}\label{Sec:EM}
\subsection{Uniformly expanding maps}\label{Sec:UE}
Several classes of uniformly expanding maps of the interval are
considered in the literature. We recall here the definition given in 
\cite{DGM} (see the references 
therein for more informations).

\begin{Definition}\label{UE}
A map $T:[0,1] \to [0,1]$ is uniformly expanding, mixing and with
density bounded from below if it satisfies the following properties:
\begin{enumerate}
\item There is a (finite or countable) partition of $T$ into
    subintervals $I_n$ on which $T$ is strictly monotonic, with a
    $C^2$ extension to its closure $\overline{I_n}$, satisfying
    Adler's condition $|T''|/|T'|^2 \leq C$, and with $|T'|\geq
    \lambda$ (where $C>0$ and $\lambda>1$ do not depend on
    $I_n$).
\item The length of $T(I_n)$ is bounded from below.
\item In this case, $T$ has finitely many absolutely continuous
    invariant measures, and each of them is mixing up to a finite
    cycle. We assume that $T$ has a single absolutely continuous
    invariant probability measure $\nu$, and that it is mixing.
\item Finally, we require that the density $h$ of $\nu$ is
    bounded from below on its support.
\end{enumerate}
\end{Definition}
From this point on, we will simply refer to such maps as
\emph{uniformly expanding}. It is well known, that, for such classes, 
the density $h$ has bounded variation.

We wish to estimates $h$ with the help of the first iterates $T, T^2, \ldots, T^n$. 
Since the bias term of a density having bounded variation is well controlled in ${\mathbb L}^1([0,1], \lambda)$, 
we shall give the rates in terms of the 
${\mathbb L}^1$-integrated risk. 
We shall use an Histogram, as defined in Section \ref{mainsec2}. 
More precisely, our estimator  $h_n$ of $h$ is given by
\begin{equation}\label{histo}
h_n(x, y)=  \sum_{i=1}^{m_n} \alpha_{i,n}(y) \varphi_{i}(x)\, , 
\end{equation}
where
$$
\varphi_i= \sqrt {m_n} {\bf 1}_{](i-1)/m_n ,  i/m_n]}  \quad  \text{and} \quad  \alpha_{i,n}(y)= \frac 1 n \sum_{k=1}^n \varphi_{i}\left (T^k (y) \right ) \, .
$$

As usual, the bias term is of order
\begin{equation} \label{biasBV}
   \int_0^1 \left |h(x)- \nu(h_n(x, \cdot)) \right | dx = O\left ( \frac 1 {m_n} \right ) \, .
\end{equation}
On another hand, one can apply Proposition \ref{p=1bis} to get
\begin{equation}\label{randomBV}
\int_0^1  \int_0^1|h_{n}(x,y)- \nu(h_n(x, \cdot))| \nu(dy) \, dx = 
O \left (\sqrt{ \frac  {m_n} {n} } \right ) \, .
\end{equation}
 Choosing $m_n=[Cn^{1/3}]$ for some $C>0$, 
 it follows from \eqref{biasBV} and \eqref{randomBV} that
$$
\int_0^1  \int_0^1|h_{n}(x,y)- h(x)| \nu(dy) \, dx = 
O \left ( n^{-1/3} \right ) \, .
$$
Now, if $\nu_{n}(y)$ is the probability measure with density $h_n(\cdot, y)$, we have just proved that
$$
\int_0^1 \|\nu_n(y)- \nu \| \,   \nu (dy) = O \left ( n^{-1/3} \right ) \, .
$$

Let us briefly explain how to derive  \eqref{randomBV} from Proposition \ref{p=1bis}. To do this, we  go back to the Markov chain associated with $T$, as we describe now. 
Let first $K$ be the
 Perron-Frobenius operator of $T$ with respect to $\nu$, defined as follows:
for any  functions $u, v$ in ${\mathbb L}^2([0,1], \nu)$
\begin{equation}\label{Perron}
\nu(u \cdot  v\circ T)=\nu(K(u) \cdot  v) \, .
\end{equation}
The relation \eqref{Perron} states  that $K$ is  the adjoint operator of the isometry $U: u \mapsto u\circ T$
acting on ${\mathbb L}^2([0,1], \nu)$. It is easy to see that the operator $K$ is a transition kernel, 
and that $\nu$ is invariant by $K$. 
Let now 
$(Y_i)_{i\geq 0}$
be a stationary Markov chain with invariant measure $\nu$
and transition kernel $K$. It is well known  that on the probability space $([0, 1], \nu)$, the
random vector $(T, T^2, \ldots , T^n)$ is distributed as
$(Y_n,Y_{n-1}, \ldots, Y_1)$. 
Hence \eqref{randomBV}  is equivalent to 
\begin{equation}\label{randomBVbis}
{\mathbb E} \left (  \int_0^1 \left | \tilde h_{n}(x)- {\mathbb E} \left (\tilde h_n(x)\right )\right |  dx \right ) = 
O \left (\sqrt{ \frac  {m_n} {n} } \right ) \, 
\end{equation}
where
$$
\tilde h_n(x)=  \sum_{i=1}^{m_n} X_{i,n} \varphi_{i}(x) \, , \quad \text{with}  \quad  X_{i,n}= \frac 1 n \sum_{k=1}^n \varphi_{i}\left (Y_k \right ) \, .
$$
Now \eqref{randomBVbis} follows easily from Proposition \ref{p=1bis} ans the fact the $\beta_{1,Y}(n)=O(a^n)$ for some $a \in (0,1)$ (see
Section 6.3 in \cite{DP2}).

\subsection{Intermittent maps}\label{Sec:LSV}
For  $\gamma$ in $(0, 1)$, we consider the intermittent map $T_\gamma$
(or simply $T$) from $[0, 1]$ to
$[0, 1]$, introduced   by Liverani, Saussol and Vaienti \cite{LSV}:
\begin{equation}\label{LSV}
   T_\gamma(x)=
  \begin{cases}
  x(1+ 2^\gamma x^\gamma) \quad  \text{ if $x \in [0, 1/2[$}\\
  2x-1 \quad \quad \quad \ \  \text{if $x \in [1/2, 1]$.}
  \end{cases}
\end{equation}
It follows from  \cite{T} that
there exists a
unique absolutely continuous $T_\gamma$-invariant probability measure $\nu_\gamma$ (or simply $\nu$),  with density $h_\gamma$ (or simply $h$). 
From \cite{T}, Theorem 1, we infer  that  the function  $x \mapsto x^\gamma h_\gamma (x)$ is bounded from above and below.
From Lemma 2.3 in \cite{LSV}, we know that $h$ is non-increasing with $h_\gamma (1)>0$, and that it is Lipshitz on any interval $[a,1]$ with $a>0$.

We wish to estimate $h$ with the help of the first iterates $T, T^2, \ldots, T^n$. To do this, we shall use the 
Histogram $h_n$ defined in \eqref{histo}.

Using the properties of $h$, it is easy to see that 
\begin{equation} \label{biasLSV}
   \int_0^1 \left |h(x)- \nu(h_n(x, \cdot)) \right | dx = O\left ( \frac 1 {m_n^{1-\gamma}} \right ) \, .
\end{equation}
On another hand, one can apply Proposition \ref{p=1bis} to get
\begin{equation}\label{randomLSV}
\int_0^1  \int_0^1|h_{n}(x,y)- \nu(h_n(x, \cdot))| \nu(dy) \, dx = 
\begin{cases}
O \left (\sqrt{ m_n/n} \right ) \qquad \ \ \,  \quad \text{ if $ \gamma < 1/2$}\\
    O \left (\sqrt{ m_n \log (n) /n} \right ) \quad \text{ if $ \gamma = 1/2$}\\
 O \left (\sqrt{ m_n/n^{(1-\gamma)/\gamma}} \right ) \quad \, \text{ if  $\gamma > 1/2$}\, .
  \end{cases}
\end{equation}
Starting from \eqref{biasLSV} and \eqref{randomLSV}, the appropriate choices of $m_n$ lead to the rates
\begin{equation}\label{ratesLSV}
\int_0^1  \int_0^1|h_{n}(x,y)- h(x)| \nu(dy) \, dx = 
\begin{cases}
O \left (n^{-(1-\gamma)/(3-2 \gamma)} \right )  \ \ \  \quad \text{ if $ \gamma < 1/2$}\\
    O \left ((n /\log (n))^{-1/4} \right ) \quad \   \text{ if $ \gamma = 1/2$}\\
 O \left (n^{-(1-\gamma)^2/\gamma (3-2 \gamma)}  \right ) \quad  \text{ if  $\gamma > 1/2$}\, .
  \end{cases}
\end{equation}

Keeping the same notations as in Section \ref{Sec:UE}, the bound \eqref{randomLSV} follows from 
Proposition \ref{p=1bis} by noting that the coefficients $\beta_{1,Y}(n)$ of the chain $(Y_i)_{i \geq 0}$
associated with $T$ satisfy $\beta_{1, Y}(n)=O(n^{-(1-\gamma)/\gamma})$ 
(see \cite{DDT}).

\section{Simulations}\label{Sec:Simu}
\subsection{Functions of an AR(1) process}
In this subsection, we first simulate the simple 
AR(1) process 
$$
X_{n+1}=\frac12 \left ( X_n + \varepsilon_{n+1} \right ) \, ,
$$
where $X_0$ is uniformly distributed over $[0,1]$,
and $(\varepsilon_i)_{i \geq 1}$ is a sequence of 
iid random variables with distribution ${\mathcal B}(1/2)$,
independent of $X_0$.  

One can check that the transition Kernel of this chain is
$$
 K(f)(x)= \frac 1 2 \left (f(x)+ f(x+1) \right ) \, ,
$$
and that the uniform distribution on $[0,1]$ is the unique
invariant distribution by $K$. Hence, the chain $(X_i)_{i \geq 0}$
is strictly stationary. 

It is well known that this chain is not $\alpha$-mixing 
in the sense of Rosenblatt \cite{R} (see for instance \cite{Br}). In fact, the
kernel $K$ is the Perron-Frobenius operator of the uniformly expanding map
$T_0$ defined in \eqref{LSV}, which is another way to see that this 
non-irreducible chain cannot be mixing in the sense of 
Rosenblatt. 

However, one can  prove that the coefficients $\beta_{2,X}$
of the chain $(X_i)_{i \geq 0}$ are such that 
$$
   \beta_{2, X}(k) \leq 2^{-k} 
$$
(see for instance 
Section 6.1 in \cite{DP2}). 

Let now $Q_{\mu, \sigma^2}$ be the inverse of the 
cumulative distribution function of the law 
${\mathcal N}(\mu, \sigma^2)$. Let then
$$
Y_i=Q_{\mu, \sigma^2}(X_i) \, .
$$
The sequence $(Y_i)_{i \geq 0}$ is also a stationary Markov chain (as 
an invertible function of a stationary Markov chain), and one can 
easily check that 
$\beta_{2,Y}(k)= \beta_{2,X}(k)$. By construction, 
$Y_i$ is ${\mathcal N}(\mu, \sigma^2)$-distributed, 
but the sequence $(Y_i)_{i \geq 0}$ is not a Gaussian 
process (otherwise it would be mixing in the sense of
Rosenblatt). 

Figure 1 shows two graphs of the kernel estimator of the density 
of $Y_i$, for $\mu=10$ and $\sigma^2=2$,
based on the simulated sample $Y_1, \ldots, Y_n$. The 
kernel $K$ is the Epanechnikov kernel (which is a Parzen
kernel of order 2, thus providing theoretically a good 
estimation when the density belongs to the the Sobolev space of order 2). 
Here we do not interfer, and let the 
software R choose an appropriate bandwidth, to see that
the default procedure delivers a correct estimation of
the density, even in this non-mixing framework. 

\begin{figure}[htb]
\centerline{
\includegraphics[width=8 cm, height=7.5 cm]{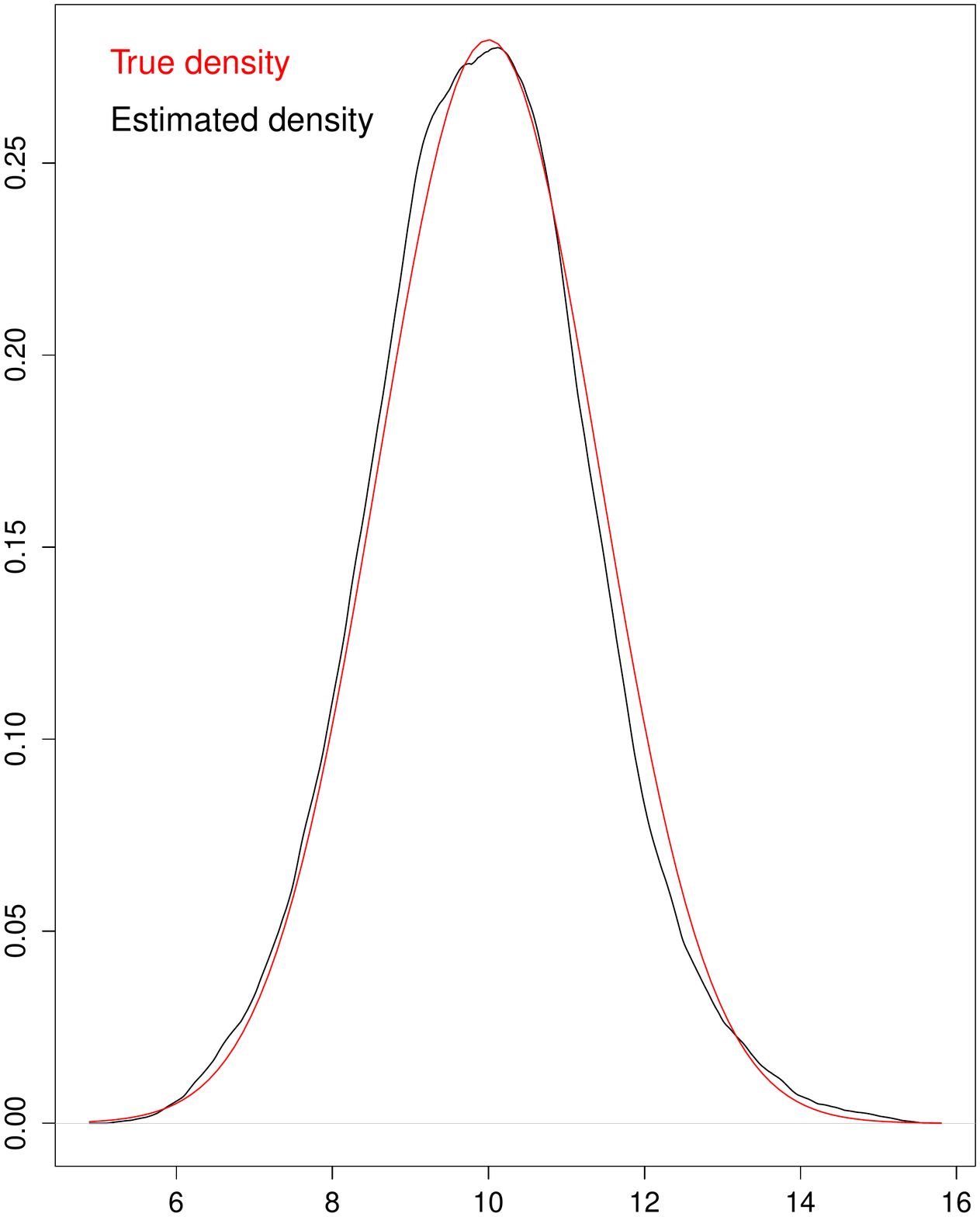}
\includegraphics[width=8 cm, height=7.5 cm]{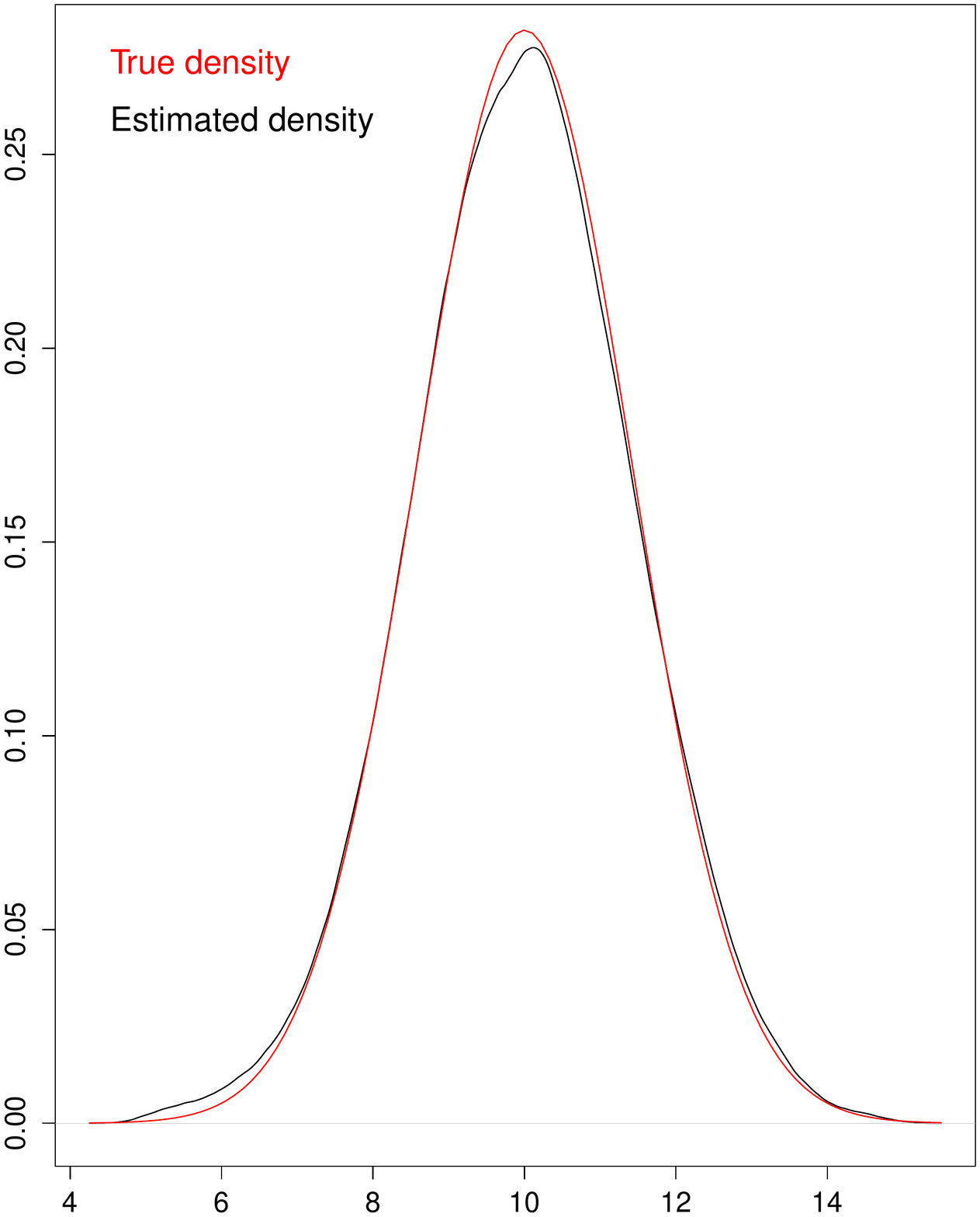}
}
\caption{Estimation of the density of $Y_i$ ($\mu=10$ and $\sigma^2=2$) $via$ the Epanechnikov kernel: $n=1000$ (left) and $n=5000$ (right)}
\label{fig:kern}
\end{figure}

We continue with another example. Let $Q: [0,1] \mapsto [0,1]$ be the inverse 
of the cumulative distribution function of the density $f$ over $[0,1]$ defined by:
$f \equiv 1/2$ on $[0,1/4] \cup [3/4,1]$ and 
$f \equiv 3/2$ on $(1/4,3/4)$. Let then
$$
Y_i=Q(X_i) \, .
$$
The same reasoning as before shows that
 $(Y_i)_{i \geq 0}$ is  a stationary Markov chain satisfying
$\beta_{2,Y}(k)= \beta_{2,X}(k)$. By construction, 
the density of the distribution of the $X_i$'s is the density $f$.
Since $f$ belongs to the class of bounded variation functions over $[0,1]$, 
the bias of the Histogram will be well controlled in ${\mathbb L}^1([0,1], \lambda)$, 
and the computations of Section \ref{mainsec2}  show that a reasonable choice 
for $m_n$ is $m_n=[n^{1/3}]$.

\begin{figure}[htb]
\centerline{
\includegraphics[width=7 cm, height=8.5 cm]{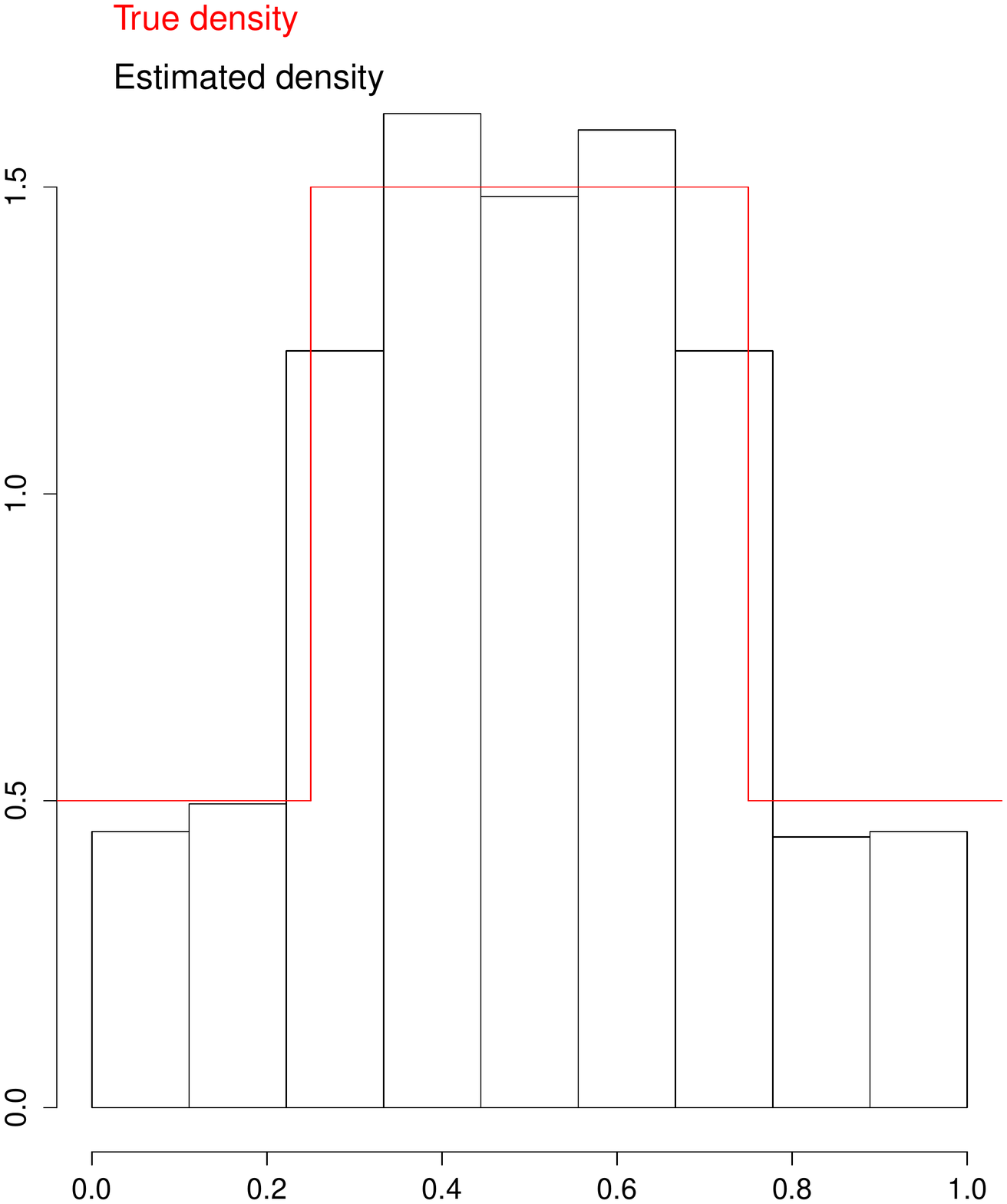}
\includegraphics[width=7 cm, height=8.5 cm]{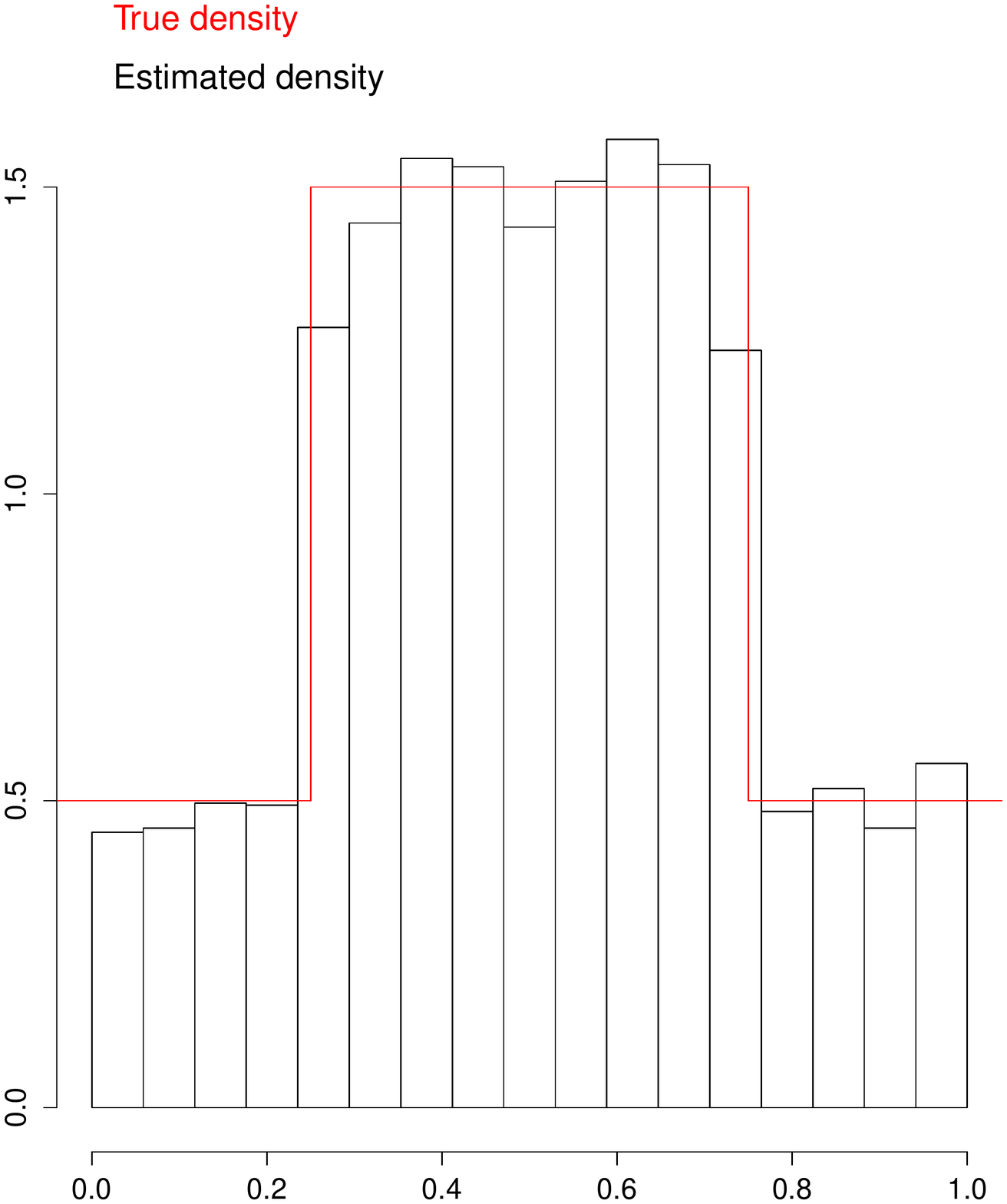}
}
\caption{Estimation of the density  $f$ by an Histogram: $n=1000$ (left) and $n=5000$ (right)}
\label{fig:hist}
\end{figure}

Figure 2 shows two Histograms
based on the simulated sample $Y_1, \ldots, Y_n$, with $m_n=[n^{1/3}]$
and two different values of $n$.
We shall now study this example from a numerical point of view, by giving an 
estimation of the ${\mathbb L}^1$-integrated risk of the Histogram. To see the 
asymptotic behavior, we let $n$ run from $5000$ to $110000$, with an increment of 
size $5000$. The ${\mathbb L}^1$-integrated risk is estimated via a classical 
Monte-Carlo procedure, by averaging the variation distance between the 
true density and the estimated density over $N=300$ independent trials. The results
are given in the table below:

\begin{center}
\begin{tabular}{|c|c||c|c|}
\hline
$n$ &    ${\mathbb L}^1$-integ. risk & $n$ & 
${\mathbb L}^1$-integ. risk \\
\hline
5000 &  0.0477   &  60000 & 0.0227\\
10000 &  0.0381  & 65000  & 0.0177\\
15000 &   0.0265& 70000 & 0.0217\\
20000 &  0.0316  & 75000 &0.0231\\
25000 &  0.0293  & 80000 & 0.0209\\
30000 &  0.0292  & 85000 &0.0202\\
35000&   0.0207 & 90000 &0.0156\\
40000 &  0.0277  & 95000 & 0.0197\\
45000 &   0.0245& 100000 &0.0209\\
50000 &  0.0191  & 105000 & 0.0193\\
55000 &  0.0251 & 110000 & 0.0189\\
\hline
\end{tabular}
\end{center}

\begin{figure}[htb]
\centerline{
\includegraphics[width=8.5 cm, height=8.3 cm]{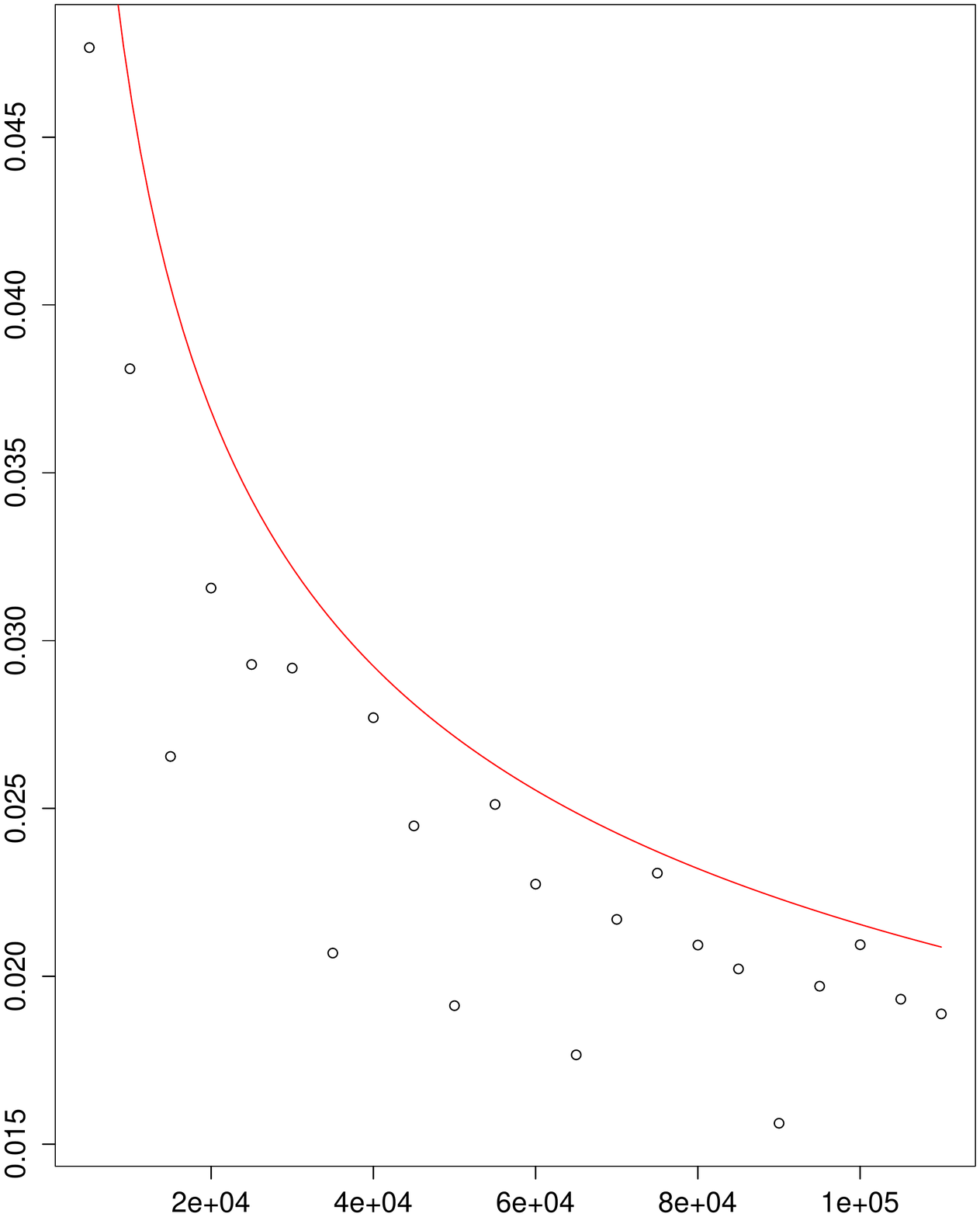}
\includegraphics[width=7 cm, height=8.3 cm]{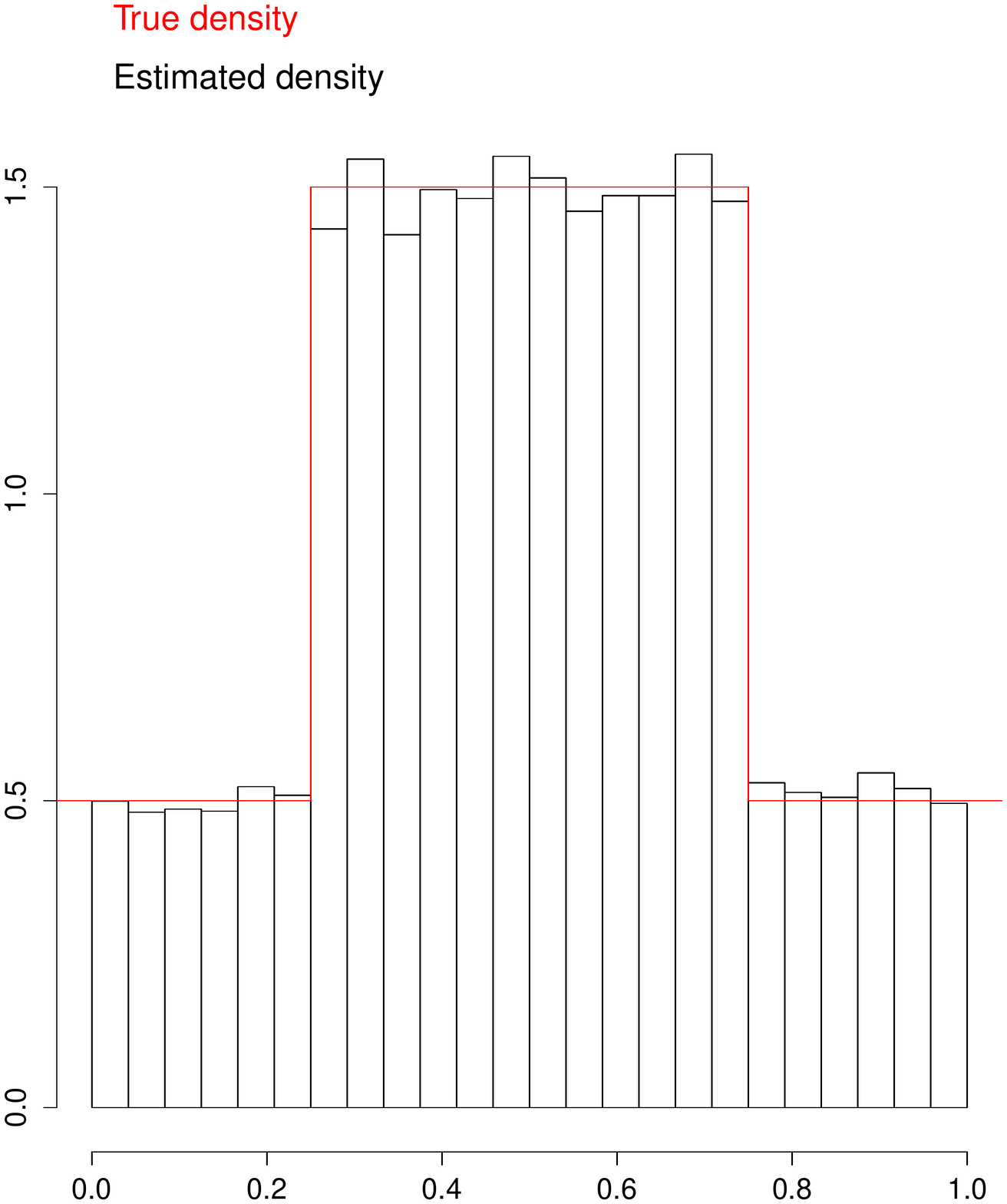}
}
\caption{Left: Graph of the ${\mathbb L}^1$-integrated risk as a function of $n$. The red curve
is the theoretical rate for estimating $BV$ functions : $n \rightarrow n^{-1/3}$.
Right: Estimation of the density  $f$ by an Histogram: $n=15000$}
\label{fig:rate}
\end{figure}

Figure 3 (left) shows the value of the 
${\mathbb L}^1$-integrated risk as $n$ increases. One can see that the ${\mathbb L}^1$-integrated risk  is smaller for some 
particular sizes of $n$. This is due to the fact that for such $n$, two of the breaks
of the Histogram are located precisely at $1/4$ and $3/4$, that is at the two discontinuity
points of the density. In that case 
the variation distance between the Histogram and the true density is particularly small, 
as illustrated by Figure 3 (right). But of course, 
we are not supposed to know where the discontinuity are located.

\subsection{Intermittent maps} \label{SimuInt}

In this section, our goal is to visualize the invariant density $h_\gamma$ of the 
intermittent maps $T_\gamma$ defined in \eqref{LSV}. Hence, we shall consider very large $n$ in order 
to get a good picture. As indicated in Subsection \ref{Sec:LSV}, we choose 
$m_n=[n^{1/(3-2\gamma)}]$ if $\gamma \in (0,1/2]$, and $m_n=[n^{(1-\gamma)/(\gamma(3-2\gamma))}]$ if $\gamma \in (1/2,1)$. 
We shall consider three cases: $\gamma=1/4$, $\gamma=1/2$ and 
$\gamma=3/4$. Recall that $\gamma <1/2$ corresponds to the short-range dependent 
case, $\gamma>1/2$ to the long-range dependent case, and $\gamma=1/2$
is the boundary case (see for instance \cite{DDT}).

\begin{figure}[htb]
\centerline{
\includegraphics[width=8.5 cm, height=8 cm]{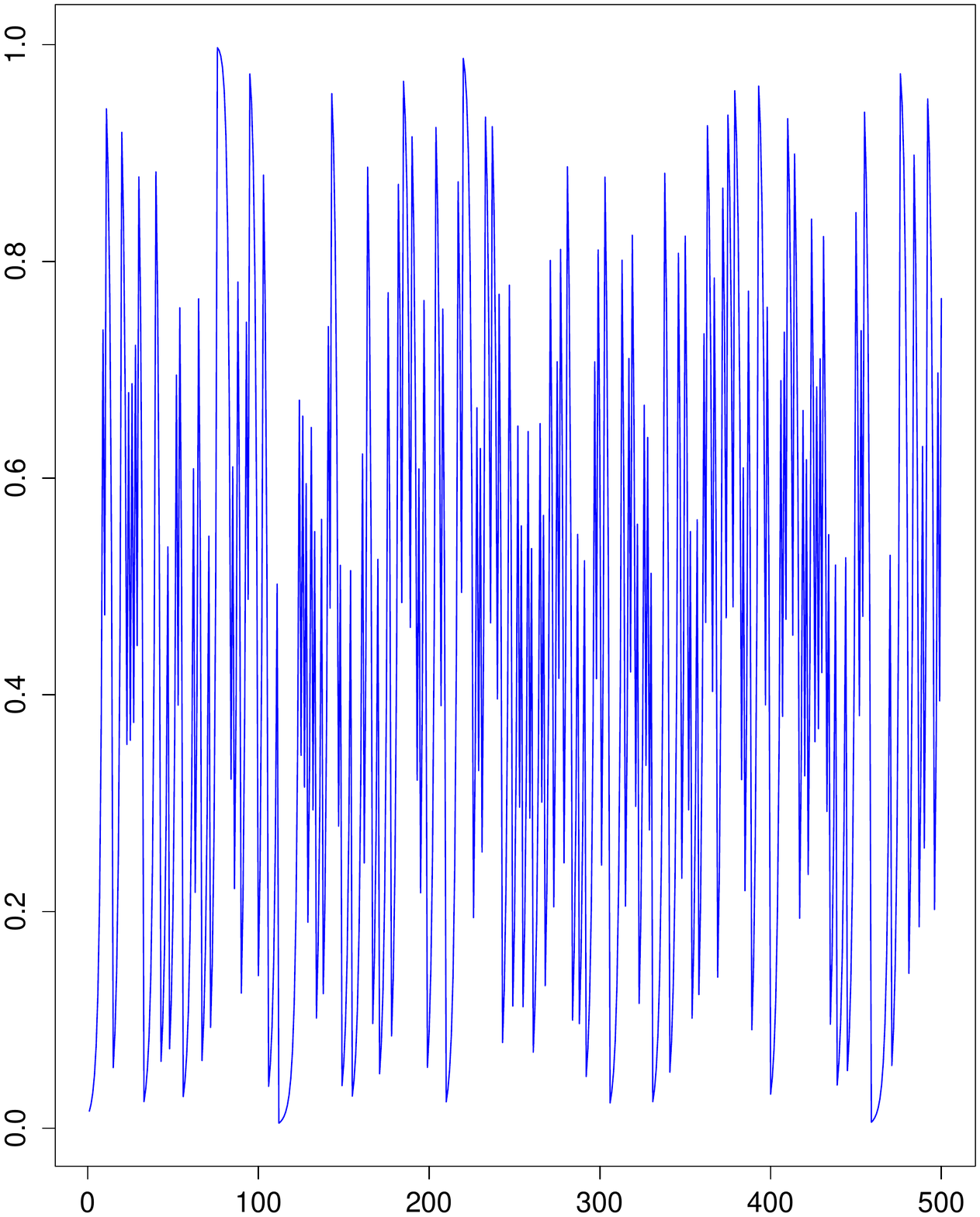}
\includegraphics[width=8.5 cm, height=8 cm]{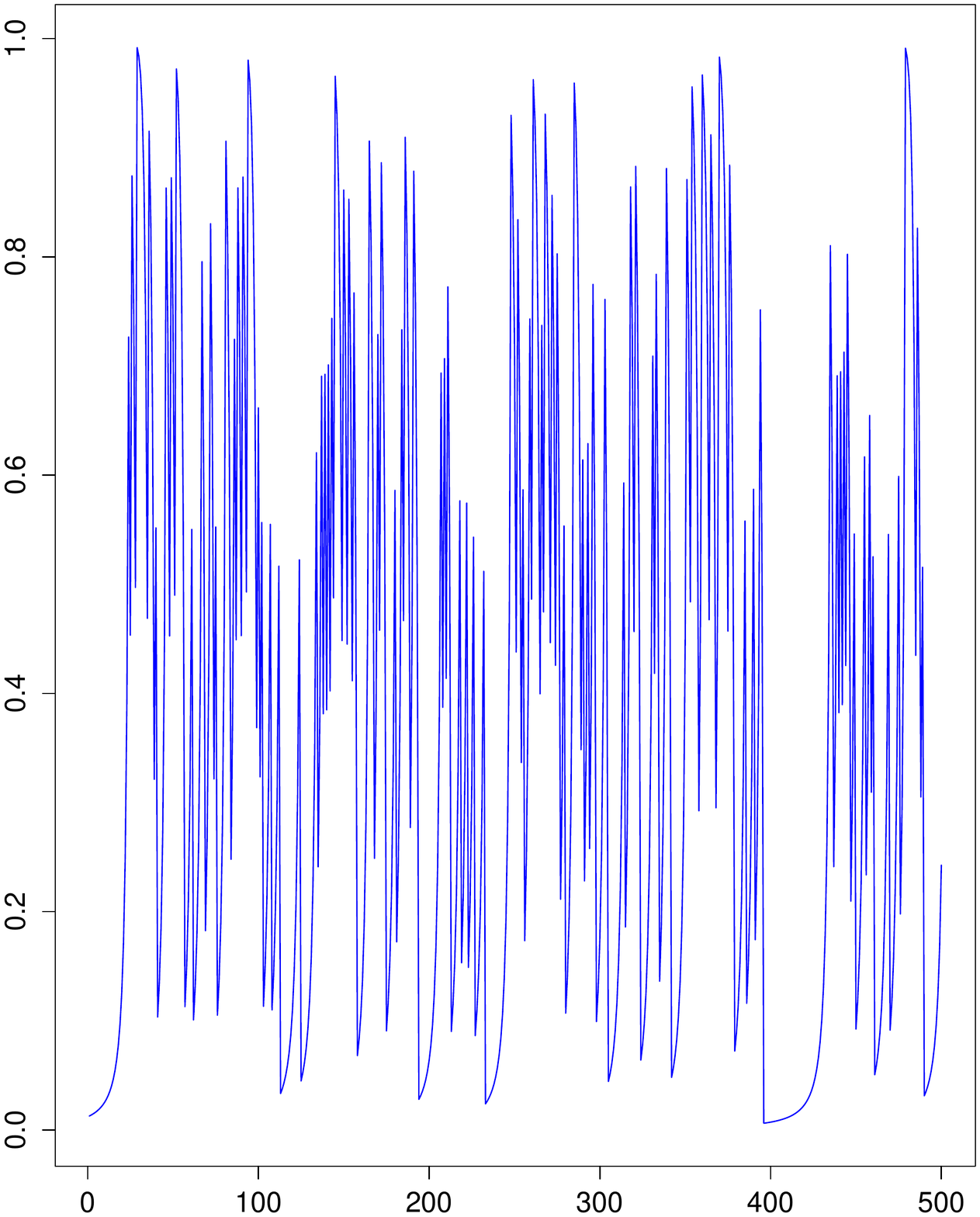}
}
\caption{Graphs of 500 iterations of the map $T_\gamma$ for $\gamma=1/4$ (left)
and $\gamma=3/4$ (right)}
\label{fig:iter}
\end{figure}

As one can see from Figure \ref{fig:iter}, due to the behavior
of $T_\gamma$ around zero, the process $(T_\gamma^i)_{i \geq 0}$ spends 
much more time in the neighborhood of $0$ when $\gamma=3/4$ than when $\gamma=1/4$. 

We do not have an explicit expression of the invariant density $h_\gamma$, but, as
already mentioned in Subsection \ref{Sec:LSV}, we know the qualitative behavior of
$h_\gamma$ in the neighborhood of $0$. More precisely, one can introduce an $equivalent$ 
density  $f_\gamma$ such that $f_\gamma(x)= (1-\gamma)x^{-\gamma}$ if $x \in (0,1]$ and $
f_\gamma(x)=0$ elsewhere.
By equivalent, we mean that there exists two positive constants $a,b,$ such that, 
on $(0,1]$,
\begin{equation}\label{equiv}
       0< a \leq h_\gamma / f_\gamma \leq b < \infty \, .
\end{equation}

\begin{figure}[htb]
\centerline{
\includegraphics[width=8.5 cm, height=6 cm]{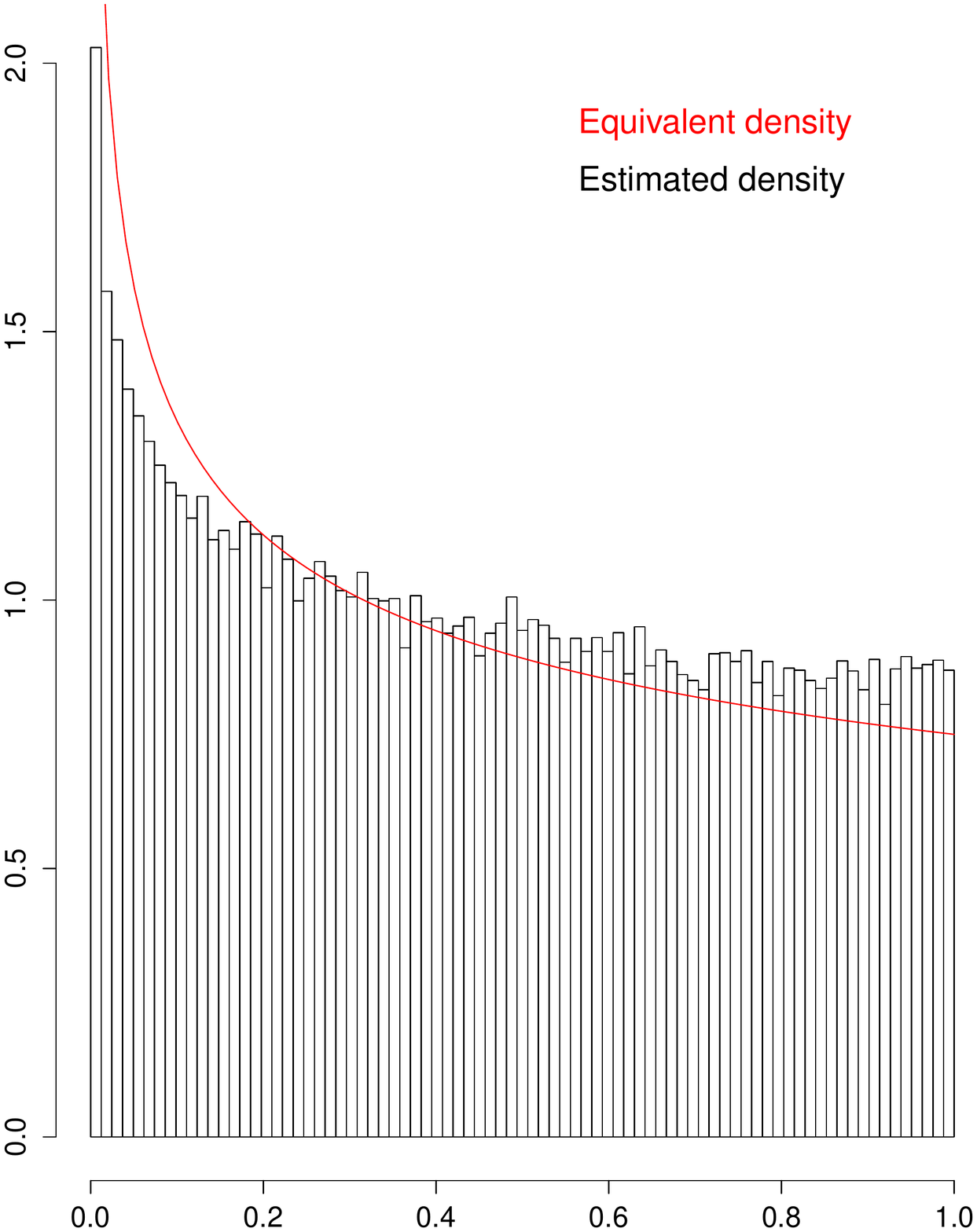}
\includegraphics[width=8.5 cm, height= 7 cm]{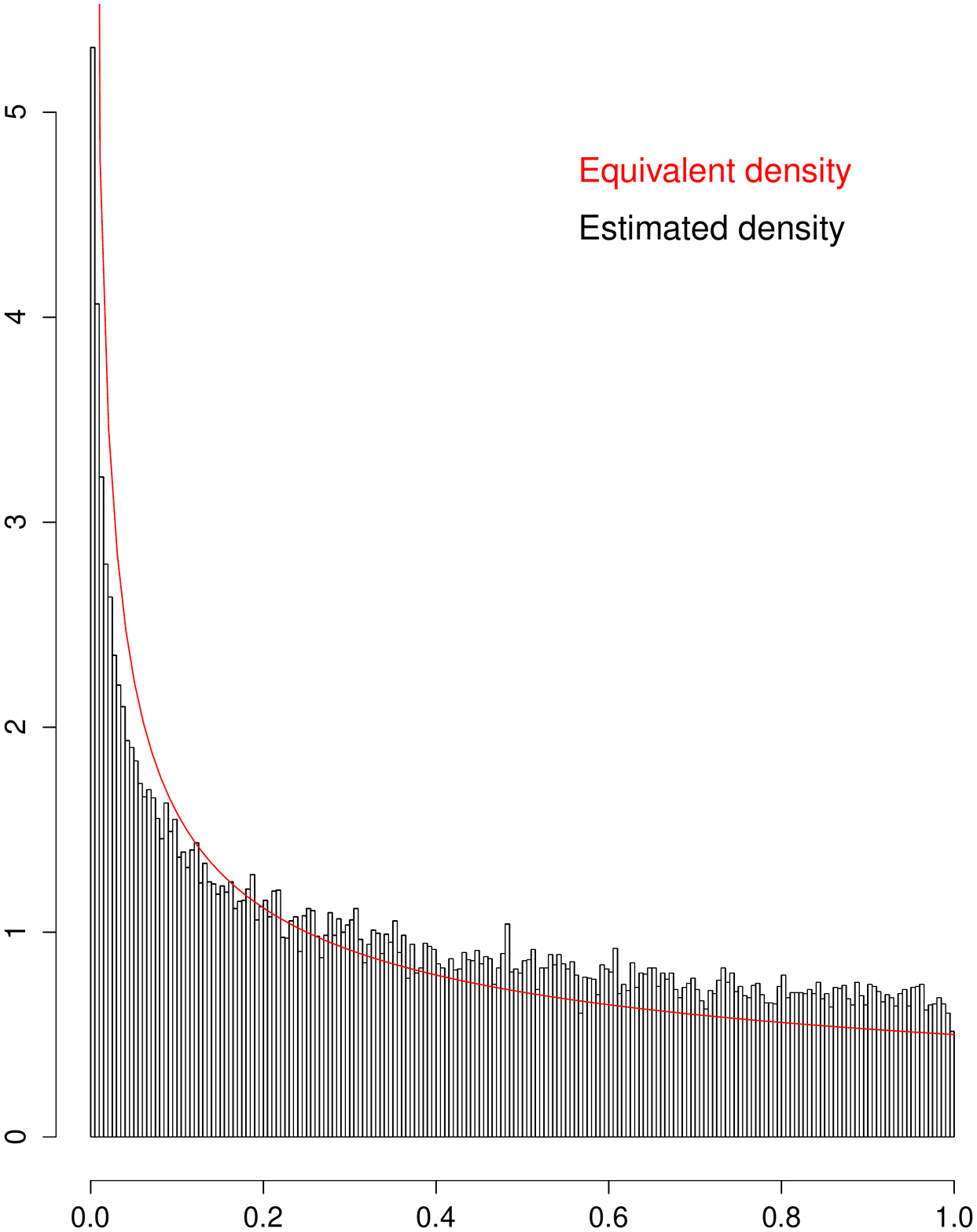}
}
\centerline{
\includegraphics[width=8.5 cm, height=8 cm]{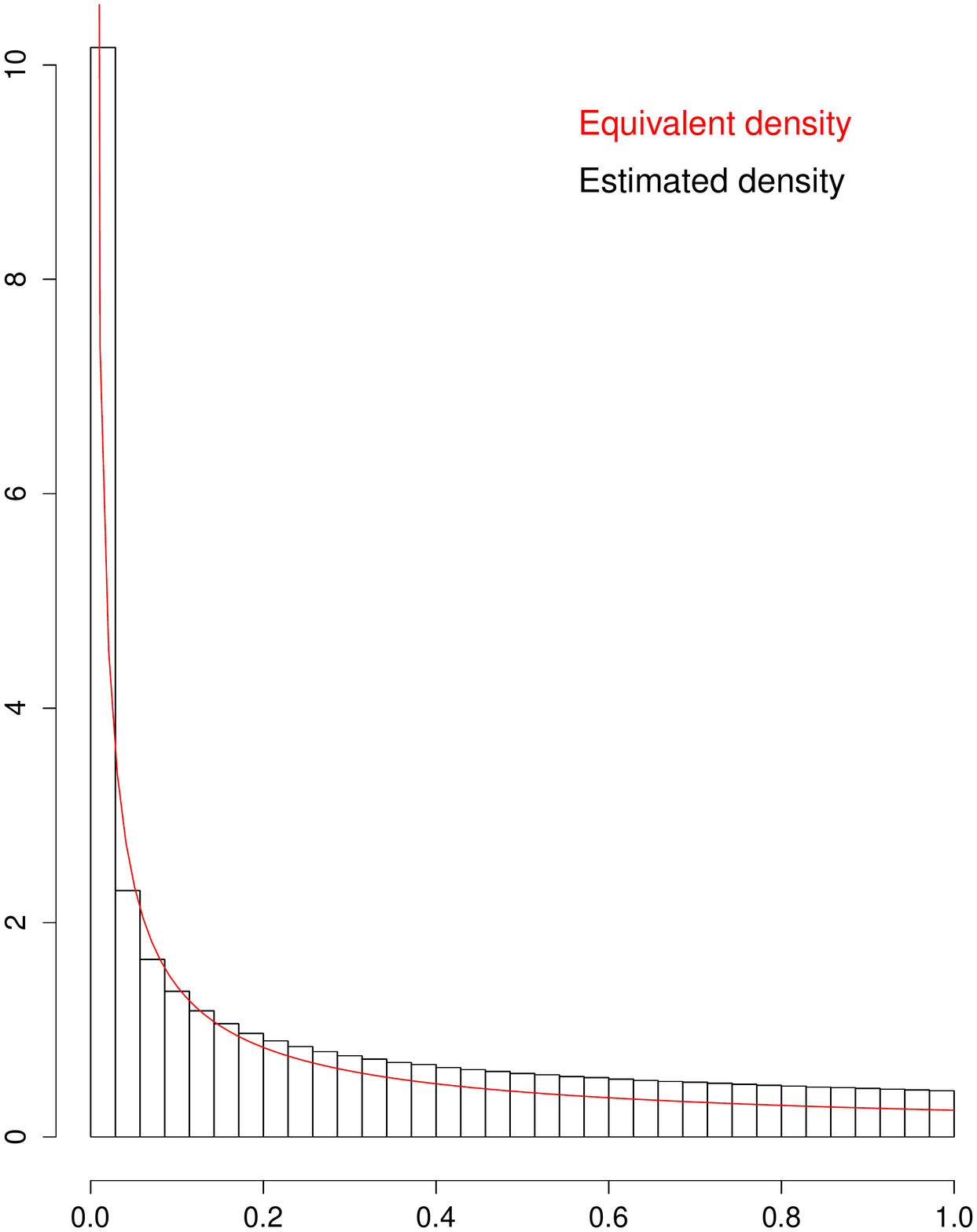}
}
\caption{Topleft: Estimation of the density $h_{1/4}$, $n=60000$.
Topright: Estimation of the density $h_{1/2}$, $n=40000$. 
Bottom: Estimation of the density $h_{3/4}$, $n=10^7$. The equivalent density is plotted in red}
\label{fig:histinter}
\end{figure}

Figure \ref{fig:histinter} shows three Histograms based on $(T_i^\gamma)_{1 \leq i \leq n}$ for 
$\gamma=1/4$, $\gamma=1/2$, $\gamma=3/4$, and very large values of $n$. Since the 
rates are very slow if $\gamma$ is much larger than $1/2$ (see \eqref{ratesLSV}), 
we have chosen $n=10^7$
for the estimation of $h_{3/4}$.
In each cases, we plotted the equivalent density on the same graph, to see 
that the behavior of $h_\gamma$ around $0$ is as expected.

\section{Proof of the results of Sections \ref{mainsec} and \ref{mainsec2} }

\subsection{Proof of Proposition \ref{p>2}}\label{proofkern1}
Setting 
$Y_{i,n}(x)=K((x-Y_{i})/h_{n})$ and
$X_{i,n}(x)=Y_{i,n}(x)-{\mathbb{E}}(Y_{i,n}(x))$, we
have that
\begin{equation}
{\mathbb{E}}\left (\int_{\mathbb{R}}|f_{n}(x)-{\mathbb{E}}(f_{n}(x))|^{p}
dx \right )\leq (nh_{n})^{-p}\int_{\mathbb{R}}{\mathbb{E}}
\left (\left  |\sum_{i=1}^{n}
X_{i,n}(x)\right  |^{p} \right ) dx\,.\label{trivialbounddensity}
\end{equation}
Starting from (\ref{trivialbounddensity}) and applying 
Proposition \ref{corRosenthalbeta}, we get
\begin{multline}\label{firstineq}
{\mathbb{E}}\left (\int_{\mathbb{R}}|f_{n}(x)-{\mathbb{E}}(f_{n}(x))|^{p}
dx \right ) \ll 
(nh_{n}^2)^{-p/2}  \int \left ( \sum_{i=0}^{n-1} \left | 
{\rm Cov} (X_{0,n}(x), X_{i,n}(x)) \right | \right )^{p/2}  dx
\\ +  n  (nh_{n})^{-p}\Vert d K \Vert^{p-1} \int \E \left (   |  Y_{0,n}(x) |  A_0(n,p)  \right ) dx
\\
+ n (nh_{n})^{-p} \Vert d K \Vert^{p-1}  \left (\int \E \left (|Y_{0,n}(x)|  \right ) dx \right ) \sum_{k=0}^{n} (k+1)^{p-2} \beta_{2,Y} (k)
\, .
\end{multline}
Since 
$$
\int |Y_{0,n}(x)| dx \leq h_n \|K\|_{1, \lambda} 
\quad \text{and} \quad 
\E(A_0(n,p)) \leq  \sum_{k=1}^{n} k^{p-2} \beta_{2,Y} (k) \, ,
$$
the two last terms on  the right hand side of  \eqref{firstineq} are bounded by the second term on the right hand side of \eqref{mainineq}.

To complete the proof, it remains to handle the first term on the right hand side of \eqref{firstineq}. By Item 1 of Lemma \ref{covbeta} (see Subsection \ref{proofrose}) applied to $Z=X_{0,n}(x)$ and 
${\mathcal F}_0= \sigma  (Y_0)$, 
\begin{multline}\label{mainkernterm}
\int \left ( \sum_{i=0}^{n-1} \left | 
{\rm Cov} (X_{0,n}(x), X_{i,n}(x)) \right | \right )^{p/2}  dx\\ 
\leq \|dK\|^{p/2}
\int \left ( \int B(y,n) |K((x-y)/h_n)| f(y) dy \right )^{p/2}  dx \, ,
\end{multline}
where $B(y,n)= b(y,0) + \cdots + b(y,n-1)$ and $b(Y_0,n) = \sup_{t \in {\mathbb R}} |P_{Y_n|Y_0}(f_t) - P(f_t)|$ (keeping the same notations as in Definition \ref{beta} for the conditional probabilities).
By Jensen's inequality 
\begin{multline*}
\left( \int B(y,n) |K((x-y)/h_n)| f(y) dy \right)^{p/2} \\ \leq h_n^{p/2} \|K\|_{1, \lambda}^{p/2-1} 
 \int B(y,n)^{p/2} f(y)^{p/2} h_n^{-1}|K((x-y)/h_n)|  dy \, .
\end{multline*}
Integrating with respect to $x$, we get 
$$
\int \left ( \int B(y,n) |K((x-y)/h_n)| f(y) dy \right )^{p/2}  dx \leq 
h_n^{p/2} \|K\|_{1, \lambda}^{p/2} 
 \int B(y,n)^{p/2} f(y)^{p/2}  dy \, .
$$
Together with \eqref{mainkernterm}, this gives 
\begin{multline}\label{main2}
(nh_{n}^2)^{-p/2}  \int \left ( \sum_{i=0}^{n-1} \left | 
{\rm Cov} (X_{0,n}(x), X_{i,n}(x)) \right | \right )^{p/2}  dx \\ \leq  (nh_{n})^{-p/2}\|dK\|^{p/2}  \|K\|_{1, \lambda}^{p/2} \int B(y,n)^{p/2} f(y)^{p/2}  dy \, .
\end{multline}
Applying H\"older's inequality,
\begin{equation}\label{end}
\int B(y,n)^{p/2} f(y)^{p/2}  dy \leq \|f \|_{p, \lambda}^{p(p-2)/2(p-1)} \left (\int  B(y,n)^{p-1} f(y) dy \right )^{p/2(p-1)}  \, .
\end{equation}
Now
$$
 \left (\int  B(y,n)^{p-1} f(y) dy\right )^{1/(p-1)}= 
 \left ({\mathbb E} \left ( B(Y_0,n)^{p-1} 
 \right)\right )^{1/(p-1)} = 
 {\mathbb E} \left ( Z_n(Y_0) B(Y_0,n) \right) \, ,
$$
where 
$$
Z_n(Y_0)=  \frac{B(Y_0,n)^{p-2}}{\left ({\mathbb E} \left ( B(Y_0,n)^{p-1} 
 \right)\right )^{(p-2)/(p-1)}} \, .
$$ 
Note that $\| Z_n(Y_0) \|_{(p-1)/(p-2)}=1$. Arguing as in \cite{Ri}  (applying Remark 1.6 with $g^2(x)=Z_n(x)$
and Inequality C.5), we infer 
that 
\begin{equation}\label{end2}
\left (\int  B(y,n)^{p-1} f(y) dy\right )^{p/2(p-1)} \ll (V_{1,p, Y}(n))^{p/2(p-1)} \, .
\end{equation}
Combining \eqref{main2}, \eqref{end} and \eqref{end2}, the proof of Proposition \ref{p>2} is complete. 

\subsection{Proof of Proposition \ref{p=1}}
We keep the same notations as in Subsection \ref{proofkern1}.
The case $p=2$ (Item 1 of Proposition \ref{p=1}) has been treated in \cite{DP}. 

For $p \in [1,2)$, we start from the elementary inequality 
$$
{\mathbb{E}}\left (\int_{\mathbb{R}}|f_{n}(x)-{\mathbb{E}}(f_{n}(x))|^{p}
dx \right )\leq
\int \left (\E \left ( |f_{n}(x)-{\mathbb{E}}(f_{n}(x))|^{2} \right ) \right )^{p/2}dx \, .
$$
Let $\alpha>1$. Applying H\"older's inequality, 
$$
{\mathbb{E}}\left (\int_{\mathbb{R}}|f_{n}(x)-{\mathbb{E}}(f_{n}(x))|^{p}
dx \right )\ll
\left (\int \left(|x|^{\alpha(2-p)/p} +1 \right) \E \left ( |f_{n}(x)-{\mathbb{E}}(f_{n}(x))|^{2} \right ) dx  \right )^{p/2} \, .
$$
As in \eqref{mainkernterm}, we infer that 
\begin{multline} \label{premier}
{\mathbb{E}}\left (\int_{\mathbb{R}}|f_{n}(x)-{\mathbb{E}}(f_{n}(x))|^{p}
dx \right )
\\
\ll  (n h_n)^{-p/2} \|dK\|^{p/2} 
\left (\int  \E \left (\left(|x|^{\alpha(2-p)/p} +1 \right) B(Y_0,n) h_n^{-1}|K((x-Y_0)/h_n)|\right ) dx  \right )^{p/2} \, .
\end{multline}
Now $ |x|^{\alpha(2-p)/p}\leq C h_n^{\alpha(2-p)/p}|(x-Y_0)/h_n|^{\alpha(2-p)/p} + C |Y_0|^{\alpha(2-p)/p}$ for some positive constant $C$. 
Plugging this upper bound in \eqref{premier} and integrating with respect to $x$ we get
\begin{multline}\label{nearend}
{\mathbb{E}}\left (\int_{\mathbb{R}}|f_{n}(x)-{\mathbb{E}}(f_{n}(x))|^{p}
dx \right )
\\
\ll
 (nh_n)^{-p/2} \|dK\|^{p/2} \|K\|_{1, \lambda}^{p/2} \left ( \E \left (|Y_0|^{\alpha(2-p)/p} B(Y_0,n) \right) \right)^{p/2} \\
+ (nh_n)^{-p/2} \|dK\|^{p/2} \left ( \left (h_n^{\alpha(2-p)/p} M_{\alpha, p}(K) + \|K\|_{1, \lambda} \right ) V_{1,2, Y}(n)\right )^{p/2}
\, .
\end{multline}
To complete the proof, it remains to handle the first term in the right hand side of \eqref{nearend}. We use once more H\"older's inequality: for any $q>1$,
\begin{multline}\label{end3}
\E \left (|Y_0|^{\alpha(2-p)/p} B(Y_0,n) \right) \leq \left(M_{\alpha q,p}(f)\right)^{1/q} 
\left ( \E\left(B(Y_0,n)^{q/(q-1)}\right)\right )^{(q-1)/q}
\\
 \ll \left(M_{\alpha q, p}(f)\right)^{1/q} \left ( U_{1, q , Y}(n)\right)^{(q-1)/q} \, ,
\end{multline}
where the last upper bound is proved as in \eqref{end2}.
Combining  \eqref{nearend}  and \eqref{end3},  the proof of Proposition \ref{p=1}
is complete. 

\subsection{Proof of Proposition \ref{p>2bis}}
\label{proofHist1}
We shall use the following notation
$$
  \nu_n(g)= \frac 1 n \sum_{k=1}^n \left ( g(Y_k)- \E(g(Y_0)) \right) \, . 
$$
With this notation, we have 
\begin{multline}
\int_0^1|f_{n}(x)-{\mathbb{E}}(f_{n}(x))|^{p} dx= \int_0^1 \sum_{j=1}^{m_n} \left | \sum_{i=1}^{r+1} \nu_n(\varphi_{i,j}) \varphi_{i,j}(x) \right |^p dx \\
\leq (r+1)^{p-1} \sum_{j=1}^{m_n} \sum_{i=1}^{r+1} \left( \int_0^1 |\varphi_{i,j}(x)|^p dx \right ) | \nu_n (\varphi_{i,j})|^p \, .
\end{multline}
Now, by definition of $\varphi_{i,j}$, 
$$
\int_0^1 |\varphi_{i,j}(x)|^p dx \leq \|R_i\|_\infty^p m_n^{p/2-1} \, .
$$
Consequently
\begin{equation}\label{step0}
{\mathbb{E}}\left (\int_0^1 |f_{n}(x)-{\mathbb{E}}(f_{n}(x))|^{p}
dx \right )\ll m_n^{p/2-1} \sum_{i=1}^{r+1}  \|R_i\|_\infty^p  \sum_{j=1}^{m_n} \E \left (| \nu_n (\varphi_{i,j})|^p \right ) \, .
\end{equation}
Applying Proposition \ref{corRosenthalbeta}, we get
\begin{multline}\label{step0bis}
\E \left (| \nu_n (\varphi_{i,j})|^p \right ) \ll 
\frac 1 {n^{p/2}} \left ( \sum_{k=0}^{n-1} \left | 
{\rm Cov} (\varphi_{i,j}(Y_0), \varphi_{i,j}(Y_k)) \right | \right )^{p/2} \\
+ \frac{1}{n^{p-1}}  \|d \varphi_{i,j} \|^{p-1} \E(\varphi_{i,j}(Y_0) A_0(n,p)) \\
+ \frac{1}{n^{p-1}}  \|d \varphi_{i,j} \|^{p-1} \E(\varphi_{i,j}(Y_0)) \sum_{k=0}^{n} (k+1)^{p-2} \beta_{2,Y} (k) \, .
\end{multline}
Since  
\begin{equation}\label{useful}
\|d \varphi_{i,j} \|= \sqrt{m_n} \|d R_i \|\, , \quad 
\sum_{j=1}^{m_n} |\varphi_{i,j}| \leq \sqrt{m_n} \|R_i\|_\infty \quad \text{and} \quad 
\E(A_0(n)) \leq \sum_{k=1}^{n} k^{p-2} \beta_{2,Y} (k) \, , 
\end{equation}
 the two last terms of the right hand side can be easily bounded, and we obtain
\begin{multline}\label{covhist}
\sum_{j=1}^{m_n}
\E \left (| \nu_n (\varphi_{i,j})|^p \right ) \ll 
\frac 1 {n^{p/2}} \sum_{j=1}^{m_n} \left ( \sum_{k=0}^{n-1} \left | 
{\rm Cov} (\varphi_{i,j}(Y_0), \varphi_{i,j}(Y_k)) \right | \right )^{p/2} \\
+ \frac{{m_n}^{p/2}}{n^{p-1}}  \|d R_i \|^{p-1} \|R_i\|_\infty \sum_{k=0}^{n} (k+1)^{p-2} \beta_{2,Y} (k) \, .
\end{multline}
It remains to control the first term on the right hand side of \eqref{covhist}. Applying Lemma \ref{covbeta} as in \eqref{mainkernterm}, we get
\begin{equation}\label{+1}
\left ( \sum_{k=0}^{n-1} \left | 
{\rm Cov} (\varphi_{i,j}(Y_0), \varphi_{i,j}(Y_k)) \right | \right )^{p/2} \leq m_n^{p/4} \|dR_i \|^{p/2} \left ( \E \left ( B(Y_0,n) \varphi_{i,j}(Y_0)\right ) \right )^{p/2}\, ,
\end{equation}
where $B(y,n)$ has been defined right after \eqref{mainkernterm}. Now, by Jensen's inequality,
\begin{equation*}
\left ( \E \left ( B(Y_0,n) \varphi_{i,j}(Y_0)\right ) \right )^{p/2}  \leq  \left ( \int |\varphi_{i,j}(x)| dx  \right)^{p/2-1} \int_0^1 B(x,n)^{p/2} f(x)^{p/2} |\varphi_{i,j}(x)| dx \, ,
\end{equation*}
and consequently, using the first part of \eqref{useful}, 
\begin{equation}\label{+2}
m_n^{p/4}\sum_{j=1}^{m_n} \left ( \E \left ( B(Y_0,n) \varphi_{i,j}(Y_0)\right ) \right )^{p/2}
\leq 
\|R_i\|_\infty^{p/2} m_n \int_0^1 B(x,n)^{p/2} f(x)^{p/2}  dx \, .
\end{equation}
From \eqref{covhist}, \eqref{+1}, \eqref{+2} and arguing as in \eqref{end}-\eqref{end2}, we get that 
\begin{multline}\label{covhistbis}
\sum_{j=1}^{m_n} \E \left (| \nu_n (\varphi_{i,j})|^p \right ) \ll 
\frac {m_n}{n^{p/2}}  \|dR_i\|^{p/2} \| R_i \|_\infty^{p/2}   \|f {\bf 1}_{[0,1]} \|_{p, \lambda}^{p(p-2)/2(p-1)} (V_{1,p, Y}(n))^{p/2(p-1)} \\
+ \frac{{m_n}^{p/2}}{n^{p-1}}  \|d R_i \|^{p-1} \|R_i\|_\infty \sum_{k=0}^{n} (k+1)^{p-2} \beta_{2,Y} (k) \, .
\end{multline}
Combining \eqref{step0} and \eqref{covhistbis}, the result follows. 

\subsection{Proof of Proposition \ref{p=1bis}}
If $p \in [1,2]$, the following inequality holds:
\begin{equation}\label{basicineq}
{\mathbb{E}}\left (\int_0^1 |f_{n}(x)-{\mathbb{E}}(f_{n}(x))|^{p}
dx \right ) \leq 
\left ({\mathbb{E}}\left (\int_0^1 |f_{n}(x)-{\mathbb{E}}(f_{n}(x))|^{2}
dx \right )\right )^{p/2} \, .
\end{equation}
To control the right hand term of \eqref{basicineq}, we apply 
\eqref{step0} with $p=2$. For $p=2$, the upper bound 
\eqref{step0bis} becomes simply
$$
\E \left (| \nu_n (\varphi_{i,j})|^2 \right ) \ll 
\frac 1 {n}  \sum_{k=0}^{n-1} \left | 
{\rm Cov} (\varphi_{i,j}(Y_0), \varphi_{i,j}(Y_k)) \right | \, ,
$$
and, following the computations of Subsection 
\ref{proofHist1}, we get 
\begin{equation}\label{covhistbis2}
\sum_{j=1}^{m_n} \E \left (| \nu_n (\varphi_{i,j})|^2 \right ) \ll 
\frac {m_n}{n}  \|dR_i\| \| R_i \|_\infty    V_{1,2, Y}(n) \, .
\end{equation} 
The result follows from \eqref{basicineq}, \eqref{step0} with $p=2$, and \eqref{covhistbis2}.

\section{Deviation and Rosenthal bounds for partial sums of 
bounded random variables }

Before proving Proposition \ref{corRosenthalbeta} in Subsection \ref{proofrose},  
we shall state and prove two 
intermediate results in Subsection \ref{devineq}:
a deviation inequality for stationary sequences of bounded random variables
(see Proposition \ref{eqmainmaxgene}) and a Rosenthal-type inequality in the 
same context (see Corollary \ref{corRosenthalborne}).

\subsection{A deviation inequality and a Rosenthal inequality} \label{devineq}
In this subsection, $(X_i)_{i \in {\mathbb Z}}$ is a 
strictly stationary sequence of real-valued
random variables such that $|X_0| \leq M$ almost surely and ${\mathbb E}(X_0)=0$. We denote by 
${\mathcal F}_i$ the $\sigma$-algebra  ${\mathcal F}_i= \sigma(X_k, k \leq i)$, 
and by ${\mathbb E}_i(\cdot)$ the conditional expectation with respect to
${\mathcal F}_i$.

\begin{Proposition} \label{deviationRosenthalgeneral}
Let $S_n = \sum_{k=1}^n X_k$. For any $x\geq M$, $r >2$, $\beta \in ]r-2,r[$ and any integer $q \in [1,n]$ such that $qM \leq x$, one has for any $x \geq M$, 
  \begin{multline}
  \label{eqmainmaxgene}
  \p \left( \sup_{ 1 \leq k \leq n}  |
  S_k | \geq 5 x \right)  \ll \frac{n^{r/2}}{x^r} \left ( \sum_{i=0}^{q-1} | {\rm Cov} (X_0, X_i)| \right )^{r/2} +   \frac{n}{ x^r } \Vert X_1 \Vert_r^r \\
  +  \frac{n}{ x^2 q }    \sum_{k=q+1}^{2q} \sum_{\ell= q +1  }^{n+q }  \left \Vert \E_{0} (X_k)   \E_{0} (X_\ell)   \right \Vert_1  \\
  + \frac{n}{x^r} q^{r/2 -1}\sum_{i=1}^{q} i^{r/2 -2}   \left \{  i \Vert  X_0  \E_{0}(X_{i}) \Vert^{r/2}_{r/2} +  \sum_{j=0}^{i-1} \Vert \E_{0}(X_i X_{i+j})  -  \E(X_i X_{i+j})\Vert_{r/2}^{r/2} \right \} \\
 +  \frac{n}{ x^r } q^{r-2 - \beta/2}\sum_{ \ell=0}^{q-1}   \sum_{j =q+1}^{n}  j^{\beta/2-1}  \Vert \E_0  ( X_j X_{j+\ell} ) - \E ( X_j X_{j + \ell} ) \Vert_{r/2}^{r/2}   \, .
  \end{multline}
\end{Proposition}
As a consequence we obtain  the following Rosenthal-type inequality:
\begin{Corollary} \label{corRosenthalborne}
 Let $S_n = \sum_{k=1}^n X_k$. For any $p > 2$, any $r \in ]2p-2,2p[$ and any 
 $\beta \in ]r-2,2p-2[$, one has
  \begin{multline}   
  \label{ineRosborne}
  \E \left( \sup_{ 1 \leq k \leq n}  |
  S_k |^p \right)  \ll   n^{p/2} \left ( \sum_{i=0}^{n-1} | {\rm Cov} (X_0, X_i)| \right )^{p/2}  \\
+  nM^{p-2} \sum_{\ell= 0  }^{n }  \sum_{k=0}^{\ell} (k+1) ^{p-3}  \left \Vert \E_{0} (X_k)   \E_{0} (X_\ell)   \right \Vert_1 \\ 
 + n M^{p-2}  \sum_{k=1}^n  k^{p-2} \left \Vert \E_{0} (X_{k})    \right \Vert_2^2 + nM^{p-r} \sum_{i=1}^n i^{p-2} \Vert X_0 \E_0(X_i) \Vert_{r/2}^{r/2}  \\
 +  n  M^{  p-r}    \sum_{j =1}^{n}  j^{\beta/2-1}  \sum_{ \ell=0}^{j-1}   ( \ell +1)^{p-2 - \beta/2} 
 \left \Vert \E_0  ( X_j X_{j+\ell} ) - \E ( X_j X_{j + \ell} ) \right \Vert_{r/2}^{r/2} \, .
  \end{multline}
\end{Corollary}

\begin{Remark}
Note that the constants that are implicitely involved in 
Proposition \ref{deviationRosenthalgeneral} 
and Corollary \ref{corRosenthalborne} depend 
only on $r$ and $\beta$. 
\end{Remark}

\noindent {\bf Proof of Proposition \ref{deviationRosenthalgeneral}}.
Let $q \in [1,n]$ be an integer such that $qM \leq x$. For any  integer $i$, define the random variables
\[ U_i =  \sum_{k=(i-1)q+1}^{iq} X_k \, .\]
Consider now the $\sigma$-algebras ${\cal G}_{i}= \F_{iq}$ and define the  variables $\tilde
U_i$ as follows:   $\tilde U_{2i -1} =  U_{2i -1}-
\E (U_{2i - 1 } | {\cal G}_{2(i-1) -1 })$  and  $ \tilde U_{2i}=U_{2i}- \E (U_{2i} | {\cal
G}_{2(i-1)}) $. The following inequality is then valid
\begin{equation*}
\max_{1 \leq k \leq  n } |S_k |\leq  2q M +
 \max_{2 \leq 2j \leq [n/q]} \left | \sum_{i=1}^j \tilde U_{2i} \right | +
\max_{1 \leq 2j-1 \leq [n/q]} \left| \sum_{i=1}^j \tilde U_{2i
-1} \right | 
  + \max_{1 \leq j \leq [n/q]} \left | \sum_{i=1}^j (U_{i} - \tilde U_{i} ) \right |  \, .
\end{equation*}
It follows that 
\begin{multline} \label{dec15FN}
\p \left ( \max_{1 \leq k \leq  n } |S_k | \geq 5x \right ) \leq  
 \p \left (\max_{2 \leq 2j \leq [n/q]} \left | \sum_{i=1}^j \tilde U_{2i} \right | \geq x \right )+
\p \left (\max_{1 \leq 2j-1 \leq [n/q]} \left | \sum_{i=1}^j \tilde U_{2i
-1} \right | \geq  x \right ) \\
  +  \p \left (\max_{1 \leq j \leq [n/q]} \Bigl | \sum_{i=1}^j (U_{i} - \tilde U_{i} ) \Bigr | \geq  x \right )  \, .
\end{multline}
Note that the two first terms on the right hand side of \eqref{dec15FN} can be treated 
similarly, so that we shall only prove an upper bound for the first one. 

Let us first deal with the last term on the right hand side of \eqref{dec15FN}. 
By Markov's inequality followed by Proposition 1 in \cite{DR}, we have
\begin{multline*}
  \p \left (\max_{1 \leq j \leq [n/q]} \left | \sum_{i=1}^j (U_{i} - \tilde U_{i} ) \right | \geq x \right ) 
    \leq \frac{4}{x^2} \sum_{i=1}^{[n/q]}  \left \Vert \E (U_{i} | {\cal
G}_{ i-2})  \right \Vert_2^2  \\
+  \frac{8}{x^2} \sum_{i=1}^{ [n/q] -1 }  \left \Vert \E (U_{i} | {\cal
G}_{i-2})  \E \left (  \sum_{j=i+1}^{[n/q]} U_{j}  | {\cal
G}_{ i-2} \right )   \right \Vert_1  \, .
\end{multline*}
Therefore, by stationarity,
\begin{multline*}
  \p \left (\max_{1 \leq j \leq [n/q]} \Bigl | \sum_{i=1}^j (U_{i} - \tilde U_{i} ) \Bigr | \geq x \right ) \\
     \leq   \frac{8}{x^2} \sum_{i=1}^{ [n/q]  }   \sum_{k=(i-1)q+1}^{iq} \sum_{j=i}^{[n/q]} \sum_{\ell=(j-1)q+1}^{jq}  \left \Vert \E_{(i-2)q} (X_k)   \E_{(i-2)q} (X_\ell)   \right \Vert_1
    \\ \leq   \frac{8}{x^2} \sum_{i=1}^{[n/q] }   \sum_{k=q+1}^{2q} \sum_{j=i}^{[n/q]} \sum_{\ell=(j-i)q + q+1  }^{(j-i +2)q}  \left \Vert \E_{0} (X_k)   \E_{0} (X_\ell)   \right \Vert_1 
     \\ \leq   \frac{8}{x^2} \sum_{i=0}^{ [n/q] -1  }   \sum_{k=q+1}^{2q} \sum_{j=0}^{i  } \sum_{\ell=( j+1) q +1  }^{(j +2)q}  \left \Vert \E_{0} (X_k)   \E_{0} (X_\ell)   \right \Vert_1  \, .
\end{multline*}
So, overall,
\begin{equation}\label{firstterm}
 \p \left (\max_{1 \leq j \leq [n/q]} \left | \sum_{i=1}^j (U_{i} - \tilde U_{i} ) \right | \geq x \right )   \leq  \frac{8 n }{ q x^2}   \sum_{k=q+1}^{2q} \sum_{\ell= q +1  }^{n+q }  \left \Vert \E_{0} (X_k)   \E_{0} (X_\ell)   \right \Vert_1 \, .
\end{equation}

Now, we handle the first term on the right hand side of \eqref{dec15FN}. Using Markov's inequality, we obtain
\beq \label{MarkovBigbloc}
 \p \left (\max_{2 \leq 2j \leq [n/q]} \left | \sum_{i=1}^j \tilde U_{2i} \right | \geq x \right ) \leq x^{-r} 
 \left \Vert\max_{2 \leq 2j \leq [n/q]} \left | \sum_{i=1}^j \tilde U_{2i} \right | \right \Vert_r^r\, .
\eeq
Note that $( \tilde U_{2i})_{i \in {\mathbb Z}} $ (resp. $( \tilde U_{2i -1
})_{i \in {\mathbb Z}}$) is a stationary sequence of martingale differences  with respect to the filtration $ ( {\cal {G}}_{2i} )_{i \in {\mathbb Z}}$ (resp. $ ( {\cal {G}}_{2i -1 } )_{i \in {\mathbb Z}}$). Applying Theorem 6 in \cite{MP}, we get
\begin{multline}\label{MP1}
\left \Vert \max_{2 \leq 2j \leq
 [n/q]
 } \left \vert \sum_{i=1}^j \tilde U_{2i}  \right \vert \right \Vert_r \\
\ll (n/q)^{1/r}  \Vert \tilde U_{2}  \Vert_r + (n/q)^{1/r} \left (  \sum_{k=1}^{[n/(2q)]} \frac{1}{k^{1+2\delta/r}}  \left \Vert  \E_0 \left ( \left ( \sum_{i=1}^k \tilde U_{2i}  \right )^2\right )\right \Vert_{r/2}^{\delta}\right )^{1/(2\delta)}\, ,
\end{multline}
where $\delta= \min ( 1, 1/(r-2))$. Since $( \tilde U_{2i})_{i \in {\mathbb Z}} $ is a stationary sequence of martingale differences  with respect to the filtration $ ( {\cal {G}}_{2i} )_{i \in {\mathbb Z}}$, 
\[
\E_0 \left ( \left ( \sum_{i=1}^k \tilde U_{2i}  \right )^2 \right ) = \sum_{i=1}^k \E_0 \left (  \tilde U^2_{2i}    \right ) \, .
\]
Moreover, 
$
\E_0  (  \tilde U^2_{2i}   ) \leq \E_0  ( U^2_{2i}    ) 
$. 
Therefore
\[
\left \Vert  \E_0 \left ( \left ( \sum_{i=1}^k \tilde U_{2i}  \right )^2\right )\right 
\Vert_{r/2} \leq \sum_{i=1}^k  \left \Vert \E_0 \left ( U^2_{2i}  \right ) - \E \left (   U^2_{2i}  \right )
\right \Vert_{r/2} +  \sum_{i=1}^k  \E \left (U^2_{2i} \right ) \, .
\]
By stationarity
\[
\sum_{i=1}^k  \E \left ( U^2_{2i}  \right ) = k \Vert S_q  \Vert_2^2  \, .
\]
Moreover
$
 \Vert \tilde U_{2}  \Vert_r \leq  2\Vert S_q \Vert_r 
$.
From \eqref{MP1} and the computations we have made, it follows that
$$ 
\left\Vert \max_{2 \leq 2j \leq
 [n/q]
 } \left \vert \sum_{i=1}^j \tilde U_{2i}  \right \vert \right \Vert^r_r \\
\ll \frac{n}{q} \Vert  S_q  \Vert^r_r + \left ( \frac{n}{q} \right )^{r/2} \Vert S_q  \Vert^r_2 +  
\frac{n}{q} \left(  \sum_{k=1}^{[n/(2q)]} \frac{1}{k^{1+2\delta/r}} D_{k,q}^{\delta}\right )^{r/(2\delta)}\, ,
$$
where 
\[
D_{k,q}  = \sum_{i=1}^k  
\left \Vert \E_0 \left ( U^2_{2i} \right ) - \E \left (   U^2_{2i}  \right ) \right \Vert_{r/2}\, .
\]
Hence,
\begin{multline}\label{consRosenthal}
\left\Vert \max_{2 \leq 2j \leq
 [n/q]
 } \left \vert \sum_{i=1}^j \tilde U_{2i}  \right \vert \right \Vert^r_r 
\ll \frac{n}{q} \Vert  S_q  \Vert^r_r +  n^{r/2}\left ( \sum_{i=0}^{q-1} | {\rm Cov} (X_0, X_i)| \right )^{r/2} 
\\+  
\frac{n}{q} \left(  \sum_{k=1}^{[n/(2q)]} \frac{1}{k^{1+2\delta/r}} D_{k,q}^{\delta}\right )^{r/(2\delta)}\, .
\end{multline}
Notice that 
\begin{multline*}
D_{k,q} \leq  \sum_{i=1}^k  \sum_{j,k =(2i-1)q+1}^{2iq}  
\left \Vert \E_0  ( X_j X_k ) - \E ( X_j X_k ) \right \Vert_{r/2} \\
 \leq 2  \sum_{i=1}^k  \sum_{j =(2i-1)q+1}^{2iq} \sum_{ \ell=0}^{2iq -j}  
 \left \Vert \E_0  ( X_j X_{j+\ell} ) - \E ( X_j X_{j + \ell} ) \right \Vert_{r/2} \, .
\end{multline*}
Let $\eta = (\beta-2)/r$ and recall that $r>2$ and
$r-2 < \beta < r$. 
Since $\eta< (r-2)/r$, applying  H\"older's inequality, we then get that
\begin{multline*}
D_{k,q}  \ll k^{- \eta +(r-2)/r } \left ( \sum_{i=1}^k i^{\beta/2-1}  
\left (  \sum_{j =(2i-1)q+1}^{2iq} \sum_{ \ell=0}^{q-1}  \Vert \E_0  ( X_j X_{j+\ell} ) - \E ( X_j X_{j + \ell} ) \Vert_{r/2} \right )^{\frac r 2}  \right )^{\frac 2 r} 
  \\
 \ll q^{2-4/r} k^{- \eta +(r-2)/r } \left ( \sum_{i=1}^k i^{\beta/2-1}    \sum_{j =(2i-1)q+1}^{2iq}  \sum_{ \ell=0}^{q-1}  \Vert \E_0  ( X_j X_{j+\ell} ) - \E ( X_j X_{j + \ell} ) \Vert_{r/2}^{r/2}  \right )^{\frac 2 r} \, .
\end{multline*}
Since $2\delta/r > - \delta \eta +\delta (r-2)/r$ (indeed $ - \eta + (r-2)/r  = (r-\beta)/r$ and $r-\beta < 2$), it follows that 
\begin{multline*}
\left ( \sum_{k=1}^{[n/(2q)]} \frac{1}{k^{1+2\delta/r}} D_{k,q}^{\delta} \right )^{r/(2 \delta)} \\
\ll  q^{r-2}  \sum_{ \ell=0}^{q-1}  \sum_{i=1}^{[n/(2q)]} i^{\beta/2-1}  \sum_{j =(2i-1)q+1}^{2iq}   \Vert \E_0  ( X_j X_{j+\ell} ) - \E ( X_j X_{j + \ell} ) \Vert_{r/2}^{r/2} 
 \\
 \ll   q^{r- \beta/2 - 1}  \sum_{ \ell=0}^{q-1}  \sum_{i=1}^{[n/(2q)]}   \sum_{j =(2i-1)q+1}^{2iq}  j^{\beta/2-1}  \Vert \E_0  ( X_j X_{j+\ell} ) - \E ( X_j X_{j + \ell} ) \Vert_{r/2}^{r/2}  \, .
\end{multline*}
So, overall, 
\beq \label{majorationMartingaleterm}
 \frac{n}{q}  \left ( \sum_{k=1}^{[n/(2q)]} \frac{1}{k^{1+2\delta/r}} D_{k,q}^{\delta} \right )^{r/(2 \delta)} 
 \ll n q^{r- \beta/2 - 2}  \sum_{ \ell=0}^{q-1}   \sum_{j =q+1}^{n}  j^{\beta/2-1}  \Vert \E_0  ( X_j X_{j+\ell} ) - \E ( X_j X_{j + \ell} ) \Vert_{r/2}^{r/2}  \, .  
 \eeq
 
 Combining \eqref{dec15FN},  \eqref{firstterm},  \eqref{consRosenthal} and 
 \eqref{majorationMartingaleterm}, Proposition \ref{deviationRosenthalgeneral} will be proved if we show that
\begin{multline} \label{majorationSqr}
 \frac{n}{q} \Vert  S_q  \Vert^r_r  \ll  n  \Vert  X_1 \Vert^r_r + n q^{r/2 -1}  \left (  \sum_{i=0}^{q-1} |\E(X_0X_i)| \right )^{r/2} \\ 
 + n q^{r/2 -1}\sum_{i=1}^{q} i^{r/2 -2}   \left \{  i \Vert  X_0  \E_{0}(X_{i}) \Vert^{r/2}_{r/2} +  \sum_{j=0}^{i-1} \Vert \E_{0}(X_i X_{i+j})  -  \E(X_i X_{i+j})\Vert_{r/2}^{r/2} \right \} \,  .
\end{multline}
By Theorem 6 in \cite{MP} again, and taking into account their Comment 7 (Item 4) together with the fact that $\Vert \E_{0}(S_{k}) \Vert_{r} \leq \Vert \E_{0}(S^2_{k}) \Vert^{1/2}_{r/2}$, we have
\[
 \Vert  S_q  \Vert^r_r \ll q  \Vert  X_1 \Vert^r_r + q \left(  \sum_{k=1}^{q}\frac{1}{k^{1+2\delta/r}}\Vert \E_{0}(S_{k}^{2}) \Vert_{r/2}^{\delta}\right)^{r/(2\delta)} \, , 
\]
where $\delta=\min(1/2,1/(p-2))$. Now,
\[
 \E(S^2_{k})  \leq 2 k \sum_{i=0}^{k-1} |\E(X_0X_i)| \, .
\]
Hence, since $r>2$, 
\beq \label{majorationSqrp1}
 \Vert  S_q  \Vert^r_r \ll q  \Vert  X_1 \Vert^r_r + q^{r/2} 
  \left (  \sum_{i=0}^{q-1} |\E(X_0X_i)| \right )^{r/2}  + q\left(  \sum_{k=1}^{q}\frac{1}{k^{1+2\delta/r}}\Vert \E_{0}(S_{k}^{2})  -  \E(S^2_{k})\Vert_{r/2}^{\delta}\right)^{r/(2\delta)} \, .
\eeq
Now
\begin{multline*}
\Vert \E_{0}(S_{k}^{2})  -  \E(S^2_{k})\Vert_{r/2} \leq 2 \sum_{i=1}^k\sum_{j=0}^{k-i} \Vert \E_{0}(X_i X_{i+j})  -  \E(X_i X_{i+j})\Vert_{r/2} \\
 \leq 2 \sum_{i=1}^k\sum_{j=0}^{i-1} \Vert \E_{0}(X_i X_{i+j})  -  \E(X_i X_{i+j})\Vert_{r/2} +   2 \sum_{i=1}^k\sum_{j=i}^{k} \Vert \E_{0}(X_i X_{i+j})  -  \E(X_i X_{i+j})\Vert_{r/2} \, .
\end{multline*}
Note that, by stationarity,
\[
\Vert \E_{0}(X_i X_{i+j})  -  \E(X_i X_{i+j})\Vert_{r/2} \leq 2 \Vert  X_i  \E_{i}(X_{i+j}) \Vert_{r/2} =  2 \Vert  X_0  \E_{0}(X_{j}) \Vert_{r/2} \, .
\]
Therefore
\[
\left \Vert \E_{0}(S_{k}^{2})  -  \E(S^2_{k}) \right \Vert_{r/2} 
\leq   2 \sum_{i=1}^k\sum_{j=0}^{i-1} \Vert \E_{0}(X_i X_{i+j})  -  \E(X_i X_{i+j})\Vert_{r/2} +   4 \sum_{j=1}^{k} j  \Vert  X_0  \E_{0}(X_{j}) \Vert_{r/2}  \, .
\]
Applying H\"older's inequality,
\[
\sum_{j=1}^{k} j  \Vert  X_0  \E_{0}(X_{j}) \Vert_{r/2} \leq   k \left ( \sum_{j=1}^{k} j^{r/2 -1}  \Vert  X_0  \E_{0}(X_{j}) \Vert^{r/2}_{r/2} \right )^{2/r}  \, , 
\]
and
\begin{multline*}
 \sum_{i=1}^k\sum_{j=0}^{i-1} \Vert \E_{0}(X_i X_{i+j})  -  \E(X_i X_{i+j})\Vert_{r/2} 
 \\ \leq k^{1-2/r} \left (  \sum_{i=1}^k \left ( \sum_{j=0}^{i-1} \Vert \E_{0}(X_i X_{i+j})  -  \E(X_i X_{i+j})\Vert_{r/2} \right )^{r/2} \right )^{2/r} \\
 \leq k^{1-2/r} \left (  \sum_{i=1}^k i^{r/2-1} \sum_{j=0}^{i-1} \Vert \E_{0}(X_i X_{i+j})  -  \E(X_i X_{i+j})\Vert_{r/2}^{r/2} \right )^{2/r}  \\
 \leq k \left (  \sum_{i=1}^k i^{r/2-2} \sum_{j=0}^{i-1} \Vert \E_{0}(X_i X_{i+j})  -  \E(X_i X_{i+j})\Vert_{r/2}^{r/2} \right )^{2/r} \, .
\end{multline*}
So, overall, 
\begin{multline*}
q\left(  \sum_{k=1}^{q}\frac{1}{k^{1+2\delta/r}}\Vert \E_{0}(S_{k}^{2})  -  \E(S^2_{k})\Vert_{r/2}^{\delta}\right)^{r/(2\delta)} \\
 \ll q^{r/2}\left \{ \sum_{j=1}^{q} j^{r/2 -1}  \Vert  X_0  \E_{0}(X_{j}) \Vert^{r/2}_{r/2} +  \sum_{i=1}^q i^{r/2-2} \sum_{j=0}^{i-1} \Vert \E_{0}(X_i X_{i+j})  -  \E(X_i X_{i+j})\Vert_{r/2}^{r/2} \right \} \, .
\end{multline*}
Taking into account this last upper bound in \eqref{majorationSqrp1}, we obtain  \eqref{majorationSqr}.   The proof of Proposition \ref{deviationRosenthalgeneral} is complete.

\medskip

\noindent{\bf Proof of Corollary \ref{corRosenthalborne}}.
Setting 
\beq \label{defofsn}
s_n^2 = \max \left (  n  \sum_{i=0}^{n-1} | {\rm Cov} (X_0, X_i)| , M^2 \right ) \, , 
\eeq
 we have 
 \begin{multline} \label{inequalityobvious}
  \E \left( \sup_{ 1 \leq k \leq n}  |
  S_k |^p \right) = p \int_0^{ nM }  x^{p-1}    \p \left( \sup_{ 1 \leq k \leq n}  |
  S_k | \geq  x \right) dx  \\ =  5^p p \int_0^{nM/5}  x^{p-1}    \p \left( \sup_{ 1 \leq k \leq n}  |
  S_k | \geq 5 x \right) dx \\
\leq    5^p p \int_0^{s_n}  x^{p-1}    \p \left( \sup_{ 1 \leq k \leq n}  |
  S_k | \geq 5 x \right) dx +  5^p p \int_{s_n}^{nM}  x^{p-1}    \p \left( \sup_{ 1 \leq k \leq n}  |
  S_k | \geq 5 x \right) dx   \, .
  \end{multline}
To handle the first term on  the right hand side of 
\eqref{inequalityobvious}, we first note that 
 \begin{multline*}
5^p p \int_0^{s_n}  x^{p-1}    \p \left( \sup_{ 1 \leq k \leq n}  |
  S_k | \geq  5 x \right) dx \leq  5^p p \int_0^{M}  x^{p-1}    \p \left( \sup_{ 1 \leq k \leq n}  |
  S_k | \geq  5 x \right) dx \\ +  5^p n^{p/2}  \left (  \sum_{i=0}^{n-1} | {\rm Cov} (X_0, X_i)|  \right )^{p/2} \, .
\end{multline*}
Now by Markov inequality followed by Proposition 1 in  \cite{DR}, we have
\begin{multline*}
5^p \int_0^{M}  x^{p-1}    \p \left( \sup_{ 1 \leq k \leq n}  |
  S_k | \geq  5 x \right) dx \\
    \leq 4 \times 5^{p-2} (p-2)^{-1} M^{p-2}  \left \{    \sum_{i=1}^n \E ( X_i^2) +  2  \sum_{i=1}^{n-1}  \left \Vert  X_i   \E_i \left (  \sum_{j=i+1}^{ n} X_j \right )   \right \Vert_1    \right \}   \, .
\end{multline*}
Hence by stationarity
\[
5^p \int_0^{M}  x^{p-1}    \p \left( \sup_{ 1 \leq k \leq n}  |
  S_k | \geq  5 x \right) dx \leq 8  \times 5^{p-2}   (p-2)^{-1} M^{p-2}   n   \sum_{j=0}^{ n -1 }  \left \Vert  X_0   \E_0  (  X_j )   \right  \Vert_1      \, .
\]
So, overall, 
 \begin{multline} \label{inequalityobvious1}
p5^p  \int_0^{s_n}  x^{p-1}    \p \left( \sup_{ 1 \leq k \leq n}  |
  S_k | \geq  x \right ) dx  \\ \leq  8p \times 5^{p-2}  n (p-2)^{-1} M^{p-2}     \sum_{j=0}^{ n -1 }  \left \Vert  X_0   \E_0  (  X_j )   \right  \Vert_1   + 5^p  n^{p/2}  \left (  \sum_{i=0}^{n-1} | {\rm Cov} (X_0, X_i)|  \right )^{p/2} \, .
  \end{multline}

We  handle now  the last term on the right hand side of \eqref{inequalityobvious}.  With this aim, we use  Proposition \ref{deviationRosenthalgeneral} with $q=[x/M]$, $r \in ]2p-2,2p[$ and  $\beta \in ]r-2,2p-2[$. 

For the first term on the right hand side of  Proposition \ref{deviationRosenthalgeneral}, we have
\begin{equation}\label{1Prop81}
n^{r/2} \int_{s_n}^{nM}  x^{p-1-r}     
\left ( \sum_{i=0}^{q-1} | {\rm Cov} (X_0, X_i)| \right )^{r/2}  dx \ll n^{p/2}  \left ( \sum_{i=0}^{n-1} | {\rm Cov} (X_0, X_i)| \right )^{p/2} \, , 
\end{equation}
since  $s_n^2 \geq   n  \sum_{i=0}^{n-1} | {\rm Cov} (X_0, X_i)| $
and $r>p$.

For the second term on the  right hand side of  Proposition \ref{deviationRosenthalgeneral}, since $r>p$ and $s_n \geq M$, 
\begin{equation}\label{2Prop81}
n \Vert X_1 \Vert_r^r  \int_{s_n}^{nM}  x^{p-1-r}      dx \ll n \Vert X_1 \Vert_r^r  s_n^{p-r} \ll n \Vert X_1 \Vert_p^p M^{r-p} s_n^{p-r} \ll n \Vert X_1 \Vert_p^p  \ll n M^{p-2}\Vert X_1 \Vert_2^2    \, .
\end{equation}

For the third term on the right hand side  of Proposition \ref{deviationRosenthalgeneral},  we have to give an upper bound for 
\[
n \int_{s_n}^{nM}  x^{p-1}    \frac{1}{ x^2 q }    \sum_{k=q+1}^{2q} \sum_{\ell= q +1  }^{n+q }  \left \Vert \E_{0} (X_k)   \E_{0} (X_\ell)   \right \Vert_1 dx \, .
\]
Write
\begin{multline*}
 \sum_{k=q+1}^{2q} \sum_{\ell= q +1  }^{n+q }  \left \Vert \E_{0} (X_k)   \E_{0} (X_\ell)   \right \Vert_1 =  \sum_{k=q+1}^{2q} \sum_{\ell= q +1  }^{n }  \left \Vert \E_{0} (X_k)   \E_{0} (X_\ell)   \right \Vert_1 \\
 +  \sum_{k=q+1}^{2q} \sum_{\ell= n +1  }^{n+q }  \left \Vert \E_{0} (X_k)   \E_{0} (X_\ell)   \right \Vert_1  \, .
\end{multline*}
Here, note that 
\[
\sum_{k=q+1}^{2q} \sum_{\ell= n +1  }^{n+q } 
 \left \Vert \E_{0} (X_k)   \E_{0} (X_\ell)   \right \Vert_1  \leq  \frac{q}{2} \sum_{k=q+1}^{2q}   
 \left \Vert \E_{0} (X_k)    \right \Vert_2^2 +  \frac{q}{2} \sum_{\ell= n +1  }^{n+q }  \left \Vert  \E_{0} (X_\ell)   \right \Vert_2^2 \, .
\]
Therefore
   \begin{multline*}
n \int_{s_n}^{nM}  x^{p-1}    \frac{1}{ x^2 q }  \sum_{k=q+1}^{2q} \sum_{\ell= n +1  }^{n+q }  \left \Vert \E_{0} (X_k)   \E_{0} (X_\ell)   \right \Vert_1 dx  \leq  \frac{n}{2} \int_{s_n}^{nM}  x^{p-3}  \sum_{k=q+1}^{2q}   \left \Vert \E_{0} (X_k)    \right \Vert_2^2  dx \\ +  \frac{n}{2} \int_{s_n}^{nM}  x^{p-3}  \sum_{\ell= n +1  }^{n+q }  \left \Vert  \E_{0} (X_\ell)   \right \Vert_2^2 dx \, .
  \end{multline*}
Now
   \begin{multline*}
n \int_{s_n}^{nM}  x^{p-3}    \sum_{k=q+1}^{2q}   \left \Vert \E_{0} (X_k)    \right \Vert_2^2  dx =  n \int_{s_n}^{nM}  x^{p-3}    \sum_{k=q+1}^{2q}   \left \Vert \E_{0} (X_k)    \right \Vert_2^2  {\bf 1}_{q \leq [n/2]}dx  \\ + n \int_{s_n}^{nM}  x^{p-3}    \sum_{k=q+1}^{2q}   \left \Vert \E_{0} (X_k)    \right  \Vert_2^2  {\bf 1}_{q  > [n/2]}dx \\
\leq   2 n \int_{s_n}^{nM}  x^{p-3}    \sum_{k=q+1}^{2q}    \left \Vert \E_{0} (X_k)    \right \Vert_2^2  {\bf 1}_{q \leq [n/2]}dx + n^2 \int_{s_n}^{nM}  x^{p-3}       \left \Vert \E_{0} (X_{[n/2 ] +1})    \right \Vert_2^2  dx \, ,
  \end{multline*}
where we  have used the fact that $q \leq n$.  Hence
 \begin{multline*}
n \int_{s_n}^{nM}  x^{p-3}    \sum_{k=q+1}^{2q}   \left \Vert \E_{0} (X_k)    \right \Vert_2^2  dx 
\leq   2 n \int_{s_n}^{nM}  x^{p-3}    \sum_{k=q+1}^{n}    \left \Vert \E_{0} (X_k)    \right \Vert_2^2  {\bf 1}_{x \leq kM}dx  \\ + n^2 \int_{s_n}^{nM}  x^{p-3}       \left  \Vert \E_{0} (X_{[n/2 ] +1})    \right \Vert_2^2  dx \\
\leq   \frac{2 n}{p-2}   M^{p-2}   \sum_{k=1}^{n}   k^{p-2} \left  \Vert \E_{0} (X_k)    \right \Vert_2^2 +   \frac{n^p}{p-2}   M^{p-2}       \left \Vert \E_{0} (X_{[n/2 ] +1})    \right \Vert_2^2    \, .
  \end{multline*}
 Now
  \[
  n^{p-1} \leq 2^{p-1}  ( [n/2] +1 )^{p-1} \leq  2^{p-1} (p-1) \sum_{k=1}^{[n/2] +1} k^{p-2} \, ,
  \]
and therefore,
\[
n^p  \left \Vert \E_{0} (X_{[n/2 ] +1})    \right \Vert_2^2 \leq 2^{p-1} (p-1) n  \sum_{k=1}^{n} k^{p-2}  
\left \Vert \E_{0} (X_{k})    \right  \Vert_2^2 \, .
\]
Consequently
\[
n \int_{s_n}^{nM}  x^{p-3}    \sum_{k=q+1}^{2q}   \left \Vert \E_{0} (X_k)    \right  \Vert_2^2  dx  \ll n M^{p-2}  \sum_{k=1}^{n}  k^{p-2} \left  \Vert \E_{0} (X_k)    \right \Vert_2^2  \, .
\]
On another hand, since $q \leq n$, 
\[
n \int_{s_n}^{nM}  x^{p-3}  \sum_{\ell= n +1  }^{n+q }  \left \Vert  \E_{0} (X_\ell)   \right \Vert_2^2 dx \leq \frac{n^p M^{p-2} }{p-2} \left \Vert \E_{0} (X_{n+1})    \right \Vert_2^2  \, .
\]
  Proceeding as before, we get
  \[
  n \int_{s_n}^{nM}  x^{p-3}  \sum_{\ell= n +1  }^{n+q }  \left \Vert  \E_{0} (X_\ell)   \right \Vert_2^2 dx \ll n M^{p-2}   \sum_{k=1}^{n}  k^{p-2} \left \Vert \E_{0} (X_k)    \right \Vert_2^2  \, . 
  \]
So, overall,
\[
n \int_{s_n}^{nM}  x^{p-1}    \frac{1}{ x^2 q }  \sum_{k=q+1}^{2q} \sum_{\ell= n +1  }^{n+q }  \big \Vert \E_{0} (X_k)   \E_{0} (X_\ell)   \big \Vert_1 dx \ll  n M^{p-2}   \sum_{k=1}^{n}  k^{p-2} \big \Vert \E_{0} (X_k)    \big \Vert_2^2  \, . 
\]  
We handle now the quantity
\[
n \int_{s_n}^{nM}  x^{p-1}    \frac{1}{ x^2 q }  \sum_{k=q+1}^{2q} \sum_{\ell= q +1  }^{n }  \left \Vert \E_{0} (X_k)   \E_{0} (X_\ell)   \right  \Vert_1 dx \, .
\]
We note first that 
 \begin{multline*}
n \int_{s_n}^{nM}  x^{p-1}    \frac{1}{ x^2 q }   \sum_{k=q+1}^{2q} \sum_{\ell= q +1  }^{n }  \left \Vert \E_{0} (X_k)   \E_{0} (X_\ell)   \right \Vert_1 dx  \\
\leq  2n \int_{s_n}^{nM}  x^{p-3}      \sum_{k=q+1}^{2q} k^{-1} \sum_{\ell= q +1  }^{n }  \left \Vert \E_{0} (X_k)   \E_{0} (X_\ell)   \right \Vert_1 dx  \, .
  \end{multline*}
Now 
   \begin{multline*}
  \sum_{k=q+1}^{2q} k^{-1} \sum_{\ell= q +1  }^{n }  \left \Vert \E_{0} (X_k)   \E_{0} (X_\ell)   \right \Vert_1 \\ =  
   \sum_{k=q+1}^{ 2q } k^{-1} \sum_{\ell= q +1  }^{ k }  \left \Vert \E_{0} (X_k)   \E_{0} (X_\ell)   \right \Vert_1 
   +    \sum_{k=q+1}^{ 2q} k^{-1} \sum_{\ell= k+1   }^{n }  \left  \Vert \E_{0} (X_k)   \E_{0} (X_\ell)   \right \Vert_1 \\
   \leq  
   \sum_{k=q+1}^{ 2q }  \sum_{\ell= q +1  }^{ k }  \ell^{-1} \left \Vert \E_{0} (X_k)   \E_{0} (X_\ell)   \right\Vert_1 
   +    \sum_{\ell= q+2   }^{ n  }   \sum_{k=q+1}^{  \ell -1}  k^{-1} \left \Vert \E_{0} (X_k)   \E_{0} (X_\ell)   
   \right \Vert_1  \, .
  \end{multline*}
Note that 
 \begin{multline*}
n \int_{s_n}^{nM}  x^{p-3}    \sum_{\ell= q+2   }^{ n  }   \sum_{k=q+1}^{  \ell -1}  k^{-1} \left \Vert \E_{0} (X_k)   \E_{0} (X_\ell)   \right  \Vert_1 dx  \\
\leq  n    \sum_{\ell=  1  }^{ n  }   \sum_{k=1}^{  \ell }  k^{-1} \left \Vert \E_{0} (X_k)   \E_{0} (X_\ell)   \right\Vert_1    \int_{s_n}^{nM}  x^{p-3}    {\bf 1}_{x \leq kM}  dx  \\
\ll   \frac{n  M^{p-2}}{p-2}  \sum_{\ell=  1  }^{ n  }   \sum_{k=1}^{  \ell }  k^{p-3} \left \Vert \E_{0} (X_k)   \E_{0} (X_\ell)   \right \Vert_1\, .
  \end{multline*}
On the other hand
 \begin{multline*}
n \int_{s_n}^{nM}  x^{p-3}     \sum_{k=q+1}^{ 2q }  \sum_{\ell= q +1  }^{ k }  \ell^{-1} \left \Vert \E_{0} (X_k)   \E_{0} (X_\ell)   \right \Vert_1  dx  \\
= n \int_{s_n}^{nM}  x^{p-3}     \sum_{k=q+1}^{ 2q }  \sum_{\ell= q +1  }^{ k }  \ell^{-1} \left \Vert \E_{0} (X_k)   \E_{0} (X_\ell)   \right  \Vert_1  {\bf 1}_{q \leq [n/2]}  dx \\ + n \int_{s_n}^{nM}  x^{p-3}     \sum_{k=q+1}^{ 2q }  \sum_{\ell= q +1  }^{ k }  \ell^{-1} \left \Vert \E_{0} (X_k)   \E_{0} (X_\ell)   \right \Vert_1 {\bf 1}_{q > [n/2]}  dx  \\
\leq  n \int_{s_n}^{nM}  x^{p-3}     \sum_{k=q+1}^{ n }  \sum_{\ell= q +1  }^{ k }  \ell^{-1} \left \Vert \E_{0} (X_k)   \E_{0} (X_\ell)   \right \Vert_1    {\bf 1}_{x \leq  \ell M}  dx  \\
 + n \int_{s_n}^{nM}  x^{p-3}     \sum_{k=q+1}^{ 2q }  \sum_{\ell= q +1  }^{ k }  \ell^{-1} \left \Vert \E_{0} (X_k)   \E_{0} (X_\ell)   \right \Vert_1 {\bf 1}_{ [n/2] < q \leq n }  dx\, .
  \end{multline*}
Proceeding as before
\begin{multline*}
 n \int_{s_n}^{nM}  x^{p-3}     \sum_{k=q+1}^{ n }  \sum_{\ell= q +1  }^{ k }  \ell^{-1} \left \Vert \E_{0} (X_k)   \E_{0} (X_\ell)   \right  \Vert_1    {\bf 1}_{x \leq  \ell M}  dx 
 \\ \ll  n  M^{p-2}  \sum_{\ell=  1  }^{ n  }   \sum_{k=1}^{  \ell }  k^{p-3} \left  \Vert \E_{0} (X_k)   \E_{0} (X_\ell)   \right  \Vert_1
\end{multline*}
and 
 \begin{multline*}
 n \int_{s_n}^{nM}  x^{p-3}     \sum_{k=q+1}^{ 2q }  \sum_{\ell= q +1  }^{ k }  \ell^{-1} \left \Vert \E_{0} (X_k)   \E_{0} (X_\ell)   \right  \Vert_1 {\bf 1}_{ [n/2] < q \leq n }  dx \\
 \ll n^p M^{p-2}  \left  \Vert \E_{0} (X_{[n/2 ] +1})    \right  \Vert_2^2 \ll  n M^{p-2}  \sum_{k=1}^n  k^{p-2}\left  \Vert \E_{0} (X_{k})    \right \Vert_2^2 \, .
  \end{multline*}
So, overall, we obtain the following upper bound
 \begin{multline}\label{3Prop81}
n \int_{s_n}^{nM}  x^{p-1}    \frac{1}{ x^2 q }    \sum_{k=q+1}^{2q} \sum_{\ell= q +1  }^{n+q }  \big \Vert \E_{0} (X_k)   \E_{0} (X_\ell)   \big \Vert_1 dx  \\ \ll  n  M^{p-2}  \sum_{\ell=  1  }^{ n  }   \sum_{k=1}^{  \ell }  k^{p-3} \big \Vert \E_{0} (X_k)   \E_{0} (X_\ell)   \big \Vert_1 + n M^{p-2}  \sum_{k=1}^n  k^{p-2} \big \Vert \E_{0} (X_{k})    \big \Vert_2^2 \, .
  \end{multline}

For the fourth term on the right hand side  of Proposition \ref{deviationRosenthalgeneral}, setting 
\[
a(i) =    i \Vert  X_0  \E_{0}(X_{i}) \Vert^{r/2}_{r/2} +  \sum_{j=0}^{i-1} \Vert \E_{0}(X_i X_{i+j})  -  \E(X_i X_{i+j})\Vert_{r/2}^{r/2} \, , 
\]
we have 
\[
n \int_{s_n}^{nM}  x^{p-1 -r}   q^{r/2 -1}\sum_{i=1}^{q} i^{r/2 -2}a(i) dx \leq  n  M^{1-r/2} \sum_{i=1}^{n } i^{r/2 -2} a(i) \int_{s_n}^{nM}  x^{p-2 -r/2}   {\bf 1}_{ x \geq i M } dx  \, .
\]
Since $r >2p-2$ and $s_n \geq M$, it follows that 
\[
n \int_{s_n}^{nM}  x^{p-1 -r}   q^{r/2 -1}\sum_{i=1}^{q} i^{r/2 -2}a(i) dx \ll  n  M^{p-r} \sum_{i=1}^{n }  i^{p-3}  a(i)    \, .
\]
We note also that  $\beta/2 > r/2 - 1 >p-2$. Therefore 
 \begin{multline*}
  \sum_{i=1}^{n }  i^{p-3}   \sum_{j=0}^{i-1} \Vert \E_{0}(X_i X_{i+j})  -  \E(X_i X_{i+j})\Vert_{r/2}^{r/2}
  \\ = \sum_{i=1}^{n }  i^{\beta/2-1 + p-2 - \beta/2}   \sum_{j=0}^{i-1} \Vert \E_{0}(X_i X_{i+j})  -  \E(X_i X_{i+j})\Vert_{r/2}^{r/2}  \\
 \leq  \sum_{i=1}^{n }  i^{\beta/2-1}   \sum_{j=0}^{i-1} (j+1)^{p-2 - \beta/2} \Vert \E_{0}(X_i X_{i+j})  -  \E(X_i X_{i+j})\Vert_{r/2}^{r/2} \, , 
 \end{multline*}
 so that 
 \begin{multline}\label{4Prop81}
 n \int_{s_n}^{nM}  x^{p-1 -r}   q^{r/2 -1}\sum_{i=1}^{q} i^{r/2 -2}a(i) dx \ll n M^{p-r} \sum_{i=1}^n 
 i^{p-2} \Vert  X_0  \E_{0}(X_{i}) \Vert^{r/2}_{r/2} \\
 + n M^{p-r}
 \sum_{i=1}^{n }  i^{\beta/2-1}   \sum_{j=0}^{i-1} (j+1)^{p-2 - \beta/2} \Vert \E_{0}(X_i X_{i+j})  -  \E(X_i X_{i+j})\Vert_{r/2}^{r/2} \, .
 \end{multline}

Finally, for the fifth term on the right hand side of of Proposition \ref{deviationRosenthalgeneral},
 \begin{multline*}
n \int_{s_n}^{nM}  x^{p-1-r}    q^{r-2 - \beta/2}\sum_{ \ell=0}^{q-1}   \sum_{j =q+1}^{n}  j^{\beta/2-1}  \Vert \E_0  ( X_j X_{j+\ell} ) - \E ( X_j X_{j + \ell} ) \Vert_{r/2}^{r/2}  dx
\\
\leq n M^{p-1-r} \int_{s_n}^{nM}      q^{p-3 - \beta/2} \sum_{ \ell=0}^{q-1}   \sum_{j =q+1}^{n}  j^{\beta/2-1}  \Vert \E_0  ( X_j X_{j+\ell} ) - \E ( X_j X_{j + \ell} ) \Vert_{r/2}^{r/2}  dx  \, ,
  \end{multline*}
since $r>p-1$ and $q<x/M$. Now, 
since $p-3 - \beta/2 < -1$ (indeed $ \beta/2 >r/2 -1$ and $r/2 >p-1$), we get 
 \begin{multline*}
n \int_{s_n}^{nM}  x^{p-1-r}    q^{r-2 - \beta/2}\sum_{ \ell=0}^{q-1}   \sum_{j =q+1}^{n}  j^{\beta/2-1}  \Vert \E_0  ( X_j X_{j+\ell} ) - \E ( X_j X_{j + \ell} ) \Vert_{r/2}^{r/2}  dx
\\
\leq n (2M)^{ \beta/2 +3 - p } M^{p-1-r} \\ \times 
\int_{s_n}^{nM}      x^{p-3 - \beta/2} \sum_{ \ell=0}^{q-1}   \sum_{j =q+1}^{n}  j^{\beta/2-1}  \Vert \E_0  ( X_j X_{j+\ell} ) - \E ( X_j X_{j + \ell} ) \Vert_{r/2}^{r/2}  dx
\\
\leq n  (2M)^{ \beta/2 +3 - p } M^{p-1-r}  \\ \times  \sum_{ \ell=0}^{n-1} 
 \sum_{j =\ell+1}^{n}  j^{\beta/2-1}  \Vert \E_0  ( X_j X_{j+\ell} ) - \E ( X_j X_{j + \ell} ) \Vert_{r/2}^{r/2} \int_{s_n}^{nM}      x^{p-3 - \beta/2}    {\bf 1}_{x \geq ( \ell +1)M} dx \, .
\end{multline*}
Hence, since $s_n \geq M$, 
\begin{multline}\label{5Prop81}
n \int_{s_n}^{nM}  x^{p-1-r}    q^{r-2 - \beta/2}\sum_{ \ell=0}^{q-1}   \sum_{j =q+1}^{n}  j^{\beta/2-1}  \Vert \E_0  ( X_j X_{j+\ell} ) - \E ( X_j X_{j + \ell} ) \Vert_{r/2}^{r/2}  dx
\\
 \ll  n  M^{  p-r}   \sum_{ \ell=0}^{n-1}   ( \ell +1)^{p-2 - \beta/2}  \sum_{j =\ell+1}^{n}  j^{\beta/2-1}  \Vert \E_0  ( X_j X_{j+\ell} ) - \E ( X_j X_{j + \ell} ) \Vert_{r/2}^{r/2}
 \\
 \ll  n  M^{  p-r}    \sum_{j =1}^{n}  j^{\beta/2-1}  \sum_{ \ell=0}^{j-1}   ( \ell +1)^{p-2 - \beta/2} \Vert \E_0  ( X_j X_{j+\ell} ) - \E ( X_j X_{j + \ell} ) \Vert_{r/2}^{r/2} \, .
  \end{multline}
  
Corollary \ref{corRosenthalborne} follows from  \eqref{inequalityobvious}, \eqref{inequalityobvious1},  and Proposition \ref{deviationRosenthalgeneral}
combined with the bounds \eqref{1Prop81},  \eqref{2Prop81},  \eqref{3Prop81},  \eqref{4Prop81} and  \eqref{5Prop81}. 
  
\subsection{Proof of Proposition \ref{corRosenthalbeta}}\label{proofrose}

The next lemma gives covariance-type inequalities in terms of 
the variables $b_\ell(i)$ and $b_\ell(i,j)$ of  Definition \ref{beta}. 
It is almost the same as Lemma 35 in \cite{MP}, 
the only difference is that in \cite{MP}, the authors used a slightly different definition 
of the variables $b_\ell(i,j)$. The proof of the version we give here can be done 
by following the proof of Lemma 35 in \cite{MP}, and is therefore omitted. 

\begin{Lemma}
\label{covbeta}Let $Z$ be a ${\mathcal{F}}_{\ell}$-measurable real-valued random
variable and let $h$ and $g$ be two BV functions (recall that
$ \Vert dh \Vert $ is the  variation norm of the measure $dh$). Let
$Z^{(0)}=Z-{\mathbb{E}}(Z)$, $h^{(0)}(Y_{i})=h(Y_{i})-{\mathbb{E}}(h(Y_{i}))$
and $g^{(0)}(Y_{j})=g(Y_{j})-{\mathbb{E}}(g(Y_{j}))$. Define the random
variables $b_\ell(i)$ and $b_\ell(i,j)$ as in Definition \ref{beta}. Then
\begin{enumerate}
\item $\left |\E\left (Z^{(0)}h^{(0)}(Y_{i}
)\right )\right |=\left |\mathrm{Cov}(Z,h(Y_{i}))\right |\leq \Vert dh \Vert \, \E \left (|Z|b_\ell(i)\right )\,.$
\item $\left |\E \left (Z^{(0)}h^{(0)}(Y_{i})g^{(0)}(Y_{j}
)\right )\right |\leq \Vert dh \Vert \Vert dg \Vert \,\E\left (|Z| b_\ell(i,j)\right )\,.$
\end{enumerate}
\end{Lemma}

We now begin the proof of Proposition \ref{corRosenthalbeta}. Note first that if $X_i=h(Y_i)-{\mathbb E}(h(Y_i))$ for some $BV$ function $h$, then
$|X_i| \leq \|dh \|$ almost surely. 
To prove Proposition \ref{corRosenthalbeta}, we apply Corollary \ref{corRosenthalborne} with 
$M=\|dh\|$, $r \in ]\max(2p-2,4),2p[$ and  $\beta \in ]r-2,2p-2[$. We have to bound up the 
second, third, fourth and fifth terms on the right hand side of \eqref{ineRosborne}. 
Let us do this in that 
order. 

To control the second term, we note that, by stationarity, 
\begin{multline*}
\E \left |   \E_0 ( X_k)  \E_0 (X_{\ell} )  \right |  = 
\E \left | \E_{- k } (X_{\ell -k } )  \E_{-k}  ( X_{  0 }   ) \right | \\
 = \E \left (  \E_{- k } (X_{\ell -k } )  \E_{-k}  ( X_{  0 }   )  \sign \left \{   \E_{- k } (X_{\ell -k } )  \E_{-k}  ( X_{  0 }   ) \right \} \right )    \\
  = \E \left (  X_{\ell -k }  \E_{-k}  ( X_{  0 }   )  \sign \left \{   \E_{- k } (X_{\ell -k } )  \E_{-k}  ( X_{  0 }   ) 
  \right \} \right )   \, .
\end{multline*}
Hence, applying Lemma \ref{covbeta}, we get  that, for any $\ell \geq k \geq 0$, 
\[
\E \left |   \E_0 ( X_k)  \E_0 (X_{\ell} )  \right |  
\leq   \Vert dh \Vert \E \left (  \left |  \E_{-k} (X_0) \right | b_{-k } ( \ell -k )  \right ) 
\leq   \Vert dh \Vert \E \left (   |  X_0 |  b_{-k} ( \ell-k)  \right )  \, .
\]
Let then
\[
T_0(n)= \sum_{\ell=1}^{n}  \sum_{k=0}^{\ell} ( k+1)^{p-3} b_{-k} ( \ell-k)   \, ,
\]
and note  that $T_0(n)$ is a positive random variable which is ${\mathcal F}_0$-measurable and such that 
\[
\E (T_0(n)) \leq  \sum_{\ell=1}^{n}  \sum_{k=0}^{\ell} ( k+1)^{p-3} \beta_{1,Y} (\ell)  \leq  C \sum_{\ell=1}^{n}   \ell^{p-2} \beta_{1,Y} (\ell)  \, ,
\]
for some positive constant $C$. Moreover
\begin{equation}\label{2}
n\|dh\|^{p-2} \sum_{\ell= 0  }^{n }  \sum_{k=0}^{\ell} (k+1) ^{p-3}  \left \Vert \E_{0} (X_k)   \E_{0} (X_\ell)   \right \Vert_1
\ll   
 n \Vert dh \Vert^{p-1} \E \left (   |  X_0 | \left (1+ T_0(n) \right ) \right ) \, .
\end{equation}

To control the third term, we note that 
$$
\left \| \E_0(X_k) \right \|_2^2 = \E \left ( \E_{-k}(X_0) X_0 \right ) \, .
$$
Hence, applying Lemma \ref{covbeta}, we get that
$$
\left \| \E_0(X_k) \right  \|_2^2  \leq \|dh\| \E \left (   |  X_0 |  b_{-k} ( 0 )  \right ) \, .
$$
Let then 
$$
U_0(n)= \sum_{k=1}^n k^{p-2} b_{-k} ( 0 ) \, ,
$$
and note that $U_0(n)$ is a positive random variable which is ${\mathcal F}_0$-measurable and such that 
$$
\E (U_0(n)) \leq  \sum_{k=1}^{n} k^{p-2} \beta_{1,Y}(k) \, .
$$
Moreover
\begin{equation}\label{3}
 n \|dh\|^{p-2}  \sum_{k=1}^n  k^{p-2} \left \Vert \E_{0} (X_{k})    \right \Vert_2^2
 \leq  n \|dh\|^{p-1} \E \left (   |  X_0 |  U_0(n) \right ) \, .
\end{equation}

To control the fourth term, let first
$Z= |X_0|^{r/2} | \E_0 (X_i) |^{r/2 -1} \sign  \{   \E_{0}(X_i)  \}$. Then
\[
\Vert X_0 \E_0(X_i) \Vert_{r/2}^{r/2} =  \E  \left ( Z X_i \right )  =  \E  \left ( (Z-\E (Z)) X_i \right )  \, .
\]
Applying Lemma \ref{covbeta}, it follows that 
\[
\Vert X_0 \E_0(X_i) \Vert_{r/2}^{r/2}  \leq  \Vert dh \Vert \E (|Z| b_0(i) )   \leq  \Vert dh \Vert^{r-1}
\E ( |X_0| b_0(i) )     \, .
\]
Let then
\[
V_0(n)= \sum_{i=1}^n i^{p-2}  b_0(i)  \, , 
\]
and note that $V_0(n)$ is a positive random variable which is ${\mathcal F}_0$-measurable and such that  $\E (V_0(n)) \leq   \sum_{i=1}^n i^{p-2} \beta_{2,Y} (i)$. Moreover
\begin{equation}\label{4}
 n\|dh\|^{p-r} \sum_{i=1}^n i^{p-2}  \Vert X_0 \E_0(X_i) \Vert_{r/2}^{r/2}  \leq  n  \Vert dh \Vert ^{p-1}
 \E \left ( |X_0|  V_0(n) \right )   \, .
\end{equation}

To control the fifth term, note that, since $r/2-1 \geq 1$,
\begin{multline*}
 \Vert  \E_0(X_i X_{j+i}) -\E (X_i X_{j+i}) \Vert_{r/2}^{r/2} \\
 =  \E  \left ( \left |   \E_0(X_i X_{j+i}) -\E (X_i X_{j+i}) \right |^{r/2-1}  \left |   \E_0(X_i X_{j+i}) -\E (X_i X_{j+i}) \right |\right ) \\
 \leq 2^{r/2-2}  \E  \left (  |   X_i X_{j+i}|^{r/2-1}  \left |   \E_0(X_i X_{j+i}) -\E (X_i X_{j+i}) \right |\right ) \\
 + 2^{r/2-2}\left |  \E (X_i X_{j+i}) \right |^{r/2-1}\E  \left (   \left |   \E_0(X_i X_{j+i}) -\E (X_i X_{j+i}) \right |\right ) \, .
\end{multline*}
Now
\begin{multline*}
\E  \left (  |   X_i X_{j+i} |^{r/2-1}  \left |   \E_0(X_i X_{j+i}) -\E (X_i X_{j+i}) \right |\right ) 
 \\ \leq \E  \left (  |   X_i  |^{r-2}  \left |   \E_0(X_i X_{j+i}) -\E (X_i X_{j+i}) \right |\right ) + 
 \E  \left (  |   X_{j+i}  |^{r-2}  \left |   \E_0(X_i X_{j+i}) -\E (X_i X_{j+i}) \right |\right ) \, .
\end{multline*}
Let $Z= \E_{-i} ( |   X_0  |^{r-2} ) \sign  \{   \E_{-i}(X_0 X_{j}) -\E (X_0 X_{j})  \}$. Notice that 
\begin{multline*}
\E  \left (  |   X_i  |^{r -2}  \left |   \E_0(X_i X_{j+i}) -\E (X_i X_{j+i}) \right |\right )  = \E  \left (  |   X_0  |^{r-2}  \left |   \E_{-i}(X_0 X_{j}) -\E (X_0 X_{j}) \right |\right ) \\ =
 \E  \left (  \E_{-i} ( |   X_0  |^{r-2})  \left |   \E_{-i}(X_0 X_{j}) -\E (X_0 X_{j}) \right |\right ) 
= \E  \left ( (Z-\E(Z)) X_0 X_{j})\right ) 
\, .
\end{multline*}
Applying Lemma \ref{covbeta}, it follows that 
\begin{multline*}
\E  \left (  |   X_i  |^{r -2}  \left |   \E_0(X_i X_{j+i}) -\E (X_i X_{j+i}) \right |\right ) \leq \Vert dh \Vert^2 \E \left (|Z|b_{-i}(0,j)\right ) \\ \leq  \Vert dh \Vert^2 \E \left (|X_0|^{r-2}b_{-i}(0,j)\right ) \leq 
 \Vert dh \Vert^{r-1} \E \left (|X_0|b_{-i}(0,j)\right )  \, .
\end{multline*}
Similarly we get 
\[
\E  \left (  |   X_{i+j}  |^{r -2}  \left |   \E_0(X_i X_{j+i}) -\E (X_i X_{j+i}) \right |\right )  \leq \|dh\|^{r-1} 
 \E \left (|X_0|b_{-i-j}(-j,0)\right )  \, .
\]
Let then
\[
W_0(n)= \sum_{i=1}^n i^{\beta/2-1} \sum_{j=0}^{i-1}  (j+1)^{p-2-\beta/2} \left ( b_{-i}(0,j) + b_{-i-j}(-j,0)\right ) \, ,
\]
and note that $W_0(n)$ is a positive random variable which is ${\mathcal F}_0$-measurable and such that 
\begin{multline*}
\E (W_0(n)) = \sum_{i=1}^n i^{\beta/2-1} \sum_{j=0}^{i-1}  (j+1)^{p-2-\beta/2} 
\left ( \E ( b_{-i}(0,j) )+ \E( b_{-i-j}(-j,0))\right ) \\
\leq 2 \sum_{i=1}^n i^{\beta/2-1} \beta_{2,Y} (i)  \sum_{j=0}^{i-1}  (j+1)^{p-2-\beta/2} \leq C  \sum_{i=1}^n i^{p-2} \beta_{2,Y} (i) \, ,
\end{multline*}
for some positive constant $C$. Moreover
\begin{multline}\label{5}
n \|dh\|^{p-r}\sum_{i=1}^n i^{\beta/2-1} \sum_{j=0}^{i-1}  (j+1)^{p-2-\beta/2}   \Vert  \E_0(X_i X_{j+i}) -\E (X_i X_{j+i}) \Vert_{r/2}^{r/2}\\
 \ll n \| dh \|^{p-1}\E \left (|X_0| W_0(n) \right ) + n \|dh\|^{p-1} \E \left (|X_0| \right )\E \left( W_0(n) \right ) \, .
\end{multline}

To conclude the proof, let 
$
B_0(n,p)= T_0(n)+ U_0(n)+ V_0(n)+ W_0(n)
$, and note that $W_0(n)$ is a positive random variable which is ${\mathcal F}_0$-measurable and such that 
$
{\mathbb E}\left(  B_0(n,p) \right ) \leq \kappa \sum_{k=1}^n k^{p-2} \beta_{2,Y} (k)
$, for some positive constant $\kappa$.  From \eqref{ineRosborne},  \eqref{2}, \eqref{3}, \eqref{4} and \eqref{5}, and since $X_0=h(Y_0)- {\mathbb E}(h(Y_0))$,
we infer that 
\begin{multline*}
  \E \left( \sup_{ 1 \leq k \leq n}  |
  S_k |^p \right)  \ll  n^{p/2} \left ( \sum_{i=0}^{n-1} | {\rm Cov} (X_0, X_i)| \right )^{p/2}  
 +  n \Vert dh \Vert^{p-1} \E \left (   |  h(Y_0) |  B_0(n,p)  \right ) \\
+ n  \Vert dh \Vert^{p-1} \E \left (|h(Y_0)|  \right )  \sum_{k=0}^{n} (k+1)^{p-2} \beta_{2,Y} (k) \, .
  \end{multline*}
Let then $A_0(n,p) =  \kappa^{-1}B_0(n,p)$, in such a way 
that $
{\mathbb E}\left(  A_0(n,p) \right ) \leq  \sum_{k=1}^n k^{p-2} \beta_{2,Y} (k)
$. The random variable $A_0(n,p)$ 
satisfies the statement of Proposition \ref{corRosenthalbeta}, and the proof is complete.


\begin{thebibliography}{9}
\bibitem{B} H. C. P. Berbee, Random walks with stationary increments and renewal theory,
Cent. Math. Tracts, Amsterdam, 1979. 
\bibitem{Br} R. C. Bradley (1986), Basic properties of strong mixing conditions, 
{\em Dependence in probability
and statistics. A survey of recent results. Oberwolfach, 1985.} E. Eberlein and M. S. Taquu editors, 
Birk\"auser, 165-192.
\bibitem{BrHi} J. H. Bramble and S. R. Hilbert (1970), Estimation of linear functionals on Sobolev spaces with application to Fourier transforms and spline interpolation, {\em SIAM J. Numer. Anal.}
 \textbf{7} 112-124.
\bibitem{BH}  J. Bretagnolle and C. Huber (1979), Estimation des densit\'es: risque minimax, {\em
Z. Wahrsch. Verw. Gebiete} \textbf{47} 119-137.
\bibitem{DDT} J. Dedecker, H. Dehling and M. S. Taqqu (2015), weak convergence of the empirical
process of intermittent maps  in ${\mathbb L}^2$ under long-range dependence, 
{\em Stoch. Dyn.} \textbf{15} 29 pp.
\bibitem{DGM} J. Dedecker, S. Gou\"ezel and F. Merlev\`ede (2012), The almost sure invariance  principle for unbounded functions of expanding maps, {\em ALEA Lat. Am. J. Probab. Math. Stat.}
 \textbf{9}  141-163. 
\bibitem{DP}  J. Dedecker and C. Prieur (2005), New dependence coefficients. Examples and applications to statistics, 
{\em Probab. Theory  Related Fields} \textbf{132}  203-236. 
\bibitem{DP2}  J. Dedecker and C. Prieur (2007), An empirical central limit theorem for dependent
sequences, 
{\em Stochastic Process.  Appl.} \textbf{117}  121-142. 
\bibitem{DR} J. Dedecker and E. Rio (2000), On the functional central limit theorem for stationary processes, 
{\em Ann. Inst. Henri Poincar\'e Probab. Stat.}
 \textbf{36} 1-34.
 \bibitem{DVL} R. A. DeVore and G. G. Lorentz, Constructive approximation, Springer-Verlag, Berlin Heidelberg New-York, 1993.
 \bibitem{LSV} C. Liverani, B.  Saussol  and S. Vaienti (1999), A probabilistic approach to intermittency, 
{\em Ergodic Theory Dynam. Systems} \textbf{19}  671-685.
\bibitem{MP}  F. Merlev\`ede and  M. Peligrad (2013),  Rosenthal-type inequalities for the maximum of partial sums of stationary processes and examples, {\em  Ann. Probab.} \textbf{41}   914-960.
\bibitem{Ri} E. Rio (2000),   Th\'eorie asymptotique des processus al\'eatoires faiblement d\'ependants,  {\em Math\'ematiques et Applications} {\bf 31}, Springer-Verlag, Berlin.
\bibitem{R} M. Rosenblatt (1956),  A central limit theorem and a strong mixing condition, {\em Proc. Nat. Acad. Sci. U. S. A.} {\bf 42} 43-47.
\bibitem{T} M. Thaler (1980), Estimates of the invariant densities of endomorphisms with indifferent fixed points.
{\em Israel J. Math.} {\bf 37}, 303--314.
\bibitem{V} G. Viennet (1997), Inequalities for absolutely regular sequences: application to density 
estimation, 
{\em Probab. Theory  Related Fields} \textbf{107}  467-492. 

\end{thebibliography}
\end{document}